\documentclass[11pt]{amsart}

\usepackage{amsmath, amsthm, amssymb, amscd, amsfonts}
\usepackage[OT2,OT1]{fontenc}
\usepackage[all]{xypic}
\usepackage{newlfont}
\usepackage[dvips]{graphics}
\usepackage[latin1]{inputenc}
\usepackage{eucal}
\usepackage{latexsym}
\usepackage[dvips, dvipsnames, usenames]{color}

\newtheorem{thm}{Theorem}[subsection]
\newtheorem{lem}[thm]{Lemma}
\newtheorem{prop}[thm]{Proposition}
\newtheorem{cor}[thm]{Corollary}
\theoremstyle{definition}
\newtheorem{rem}[thm]{Remark}
\newtheorem{defn}[thm]{Definition}
\newtheorem{exa}[thm]{Example}

\def\bp{\begin{proof}}
\def\ep{\end{proof}}

\makeatletter
\renewcommand{\subsection}{\@startsection{subsection}{2}{0pt}{-3ex
plus -1ex minus -0.2ex}{-2mm plus -0pt minus
-2pt}{\normalfont\bfseries}}
\renewcommand{\subsubsection}{\@startsection{subsubsection}{2}{0pt}{-3ex
plus -1ex minus -0.2ex}{-2mm plus -0pt minus
-2pt}{\normalfont\bfseries}} \makeatother

\numberwithin{equation}{section}

\newcommand{\erem}{\hfill$\lozenge$\end{rem}\vskip 3pt }

\newcommand{\ord}{\operatorname{ord}}

\newcommand{\Aut}{\operatorname{Aut}}

\newcommand{\en}{\enspace }

\newcommand{\wH}{\widetilde{H}}

\newcommand{\vi}{${\en\mathsf {(i)}}\;$}
\newcommand{\vii}{${\;\mathsf {(ii)}}\;$}
\newcommand{\viii}{${\mathsf {(iii)}}\;$}

\newcommand{\Ad}{\operatorname{Ad}}

\def\C{\mathbb{C}}
\def\k{\mathbf{k}}

\def\id{\mathrm{id}}

\def\N{\mathbb{N}}

\def\dim{\mathrm{dim}}

\def\Rep{{\operatorname{Rep}}}

\def\Rad{{\operatorname{Rad}}}
\def\Fun{{\operatorname{Fun}}}

\def\Z{{\mathbb Z}}
\def\hyd{^H_H\mathcal{YD}}

\def\gyd{^{\k[\Gamma]}_{\k[\Gamma]}\mathcal{YD}}
\def\zmyd{^{\k \Z_m}_{\k \Z_m}\mathcal{YD}}
\def\uno{\mathbf{1}}

\def\bB{\mathcal{B}}
\def\bC{\mathcal{C}}
\def\bD{\mathcal{D}}
\def\bX{\mathcal{X}}
\def\bF{\mathcal{F}}
\def\bZ{\mathcal{Z}}
\def\hga{\hat{\Gamma}}
\def\aut{\underline{\operatorname{Aut}}}

\begin{document}
\title{\en Basic quasi-Hopf algebras over cyclic groups}
\author{Iv\'an Ezequiel Angiono}
\address{Facultad of Matem\'atica, Astronom\'\i a y F\'\i sica
\newline \indent
Universidad Nacional of C\'ordoba
\newline
\indent CIEM -- CONICET
\newline
\indent (5000) Ciudad Universitaria, C\'ordoba, Argentina}
 \email{angiono@mate.uncor.edu}
\date{\today}

\begin{abstract} Let $m$ a positive integer, not divisible by 2,3,5,7. We generalize the classification of basic quasi-Hopf algebras over cyclic groups of prime order given in \cite{EG3}  to the case of cyclic groups of order $m$. To this end, we introduce a family of non-semisimple radically graded quasi-Hopf algebras $A(H,s)$, constructed as subalgebras of Hopf algebras twisted by a quasi-Hopf twist, which are not twist equivalent to Hopf algebras. Any basic quasi-Hopf algebra over a cyclic group of order $m$ is either semisimple, or is twist equivalent to a Hopf algebra or a quasi-Hopf algebra of type $A(H,s)$.
\end{abstract}

\maketitle

\section{Introduction}\label{section:introduction}

A finite dimensional associative algebra is \emph{basic} if
all its irreducible representations are
1-dimensional. Dually, we obtain pointed coalgebras. Thus, the problem classification
of basic Hopf algebras up to isomorphism is equivalent to
the problem of classification of finite dimensional pointed Hopf algebras up to isomorphism.
When the group $G(H)$ of grouplike elements of a
finite dimensional pointed Hopf algebra $H$ is abelian of order not divisible by $2,3,5,7$,
this problem was solved by Andruskiewitsch and Schneider, see \cite{AS4}.
One of the main difficulties is, once one knows all the coradically
graded pointed Hopf algebras which are finite dimensional, to obtain
all the liftings; i.e. for any coradically graded $H_0$, to find all
the Hopf algebras $H$ whose associated graded Hopf algebra
is $H_0$.

The result of Andruskiewitsch-Schneider also yields a classification of pointed finite tensor
categories with abelian groups of grouplike elements of order not divisible by $2,3,5,7$
which have a fiber functor, as the categories of comodules over such
pointed Hopf algebras (see \cite{EO}). Moreover, by Masouka's Theorem \cite[Thm. A1]{Ma}, the equivalence classes of
such categories reduce to the graded case, because the category of comodules
over a lifting $H$ of $H_0$ is equivalent to the category of comodules over $H_0$.

In what follows, $\k$ will denote an algebraically closed field of characteristic zero. All the algebras and tensor categories considered in this work are over $\k$.

The general problem of classification of pointed finite tensor categories (not necessarily having a fiber functor) reduces to classification of basic quasi-Hopf algebras up to twist, and it is closely related to the classification of pointed Hopf algebras. The first approach to this problem was suggested by Etingof and Gelaki in a series of papers, in which they classified pointed finite tensor categories whose group of invertible objects has prime order (see \cite{EG3} for a complete answer). To do so, they considered the quasi-Hopf algebras $A(q)$ with non-trivial associator constructed in \cite{G}. This family completes the list of such categories, with the categories of representations of Hopf algebras and the semisimple ones.

In this work we classify pointed tensor categories such that their group of invertible objects is cyclic, with order not divisible by $2,3,5,7$. This restriction on the order comes mainly from the classification Theorem for pointed Hopf algebras over abelian groups in \cite{AS4}. In fact, we can classify basic radically graded quasi-Hopf algebras over cyclic groups of odd order, but restrict as above when we consider liftings of these algebras. Therefore main Theorem could still hold for any cyclic group of odd order if one can extend the theory of liftings for any group of odd order.

This family of categories has a subfamily corresponding to non-semisimple quasi-Hopf algebras $A(H,s)$, constructed in a similar way to the family of quasi-Hopf algebras $A(q)$ of Gelaki, from radically graded Hopf algebras. Consider $H=\oplus_{n \geq 0} H(n)$ a radically graded Hopf algebra, generated by a group like element $\chi$ of order $m^{2}$ and skew primitive elements $x_1,...,x_{\theta}$ such that:
\begin{equation}\label{skewprimitives}
\chi x_i\chi^{-1} = q^{d_i}x_i, \quad \Delta(x_i)= x_i \otimes \chi^{b_i} + 1 \otimes x_i,
\end{equation}
where $q$ is a root of unity of order $m^{2}$. Call
\begin{equation}\label{solutionsHforpossiblecocycles}
    \Upsilon(H):= \left\{ s \in \{1, \ldots, m-1\}: b_i \equiv sd_i (m), \, 1 \leq i \leq \theta \right\}.
\end{equation}
Consider its subalgebra $A(H,s)$ generated by $\sigma:=\chi^{m}$ and $x_1,...,x_{\theta}$. Modifying the coalgebra structure of $H$ by a twist $J_s \in H\otimes H$ (there exists one $J_s$ for each $s\in \Upsilon(H)$), we shall prove that $A(H,s)$ with the induced coalgebra structure by restriction is a quasi Hopf algebra, which is not twist equivalent to a Hopf algebra.

As we will prove that liftings of quasi-Hopf algebras $A(H,s)$ come from de-equivariantizations of liftings of Hopf algebras, Masuoka's Theorem simplifies the classification problem: we can restrict to the radically graded case.

The main result of this work is the following:
\begin{thm}\label{thm:classificationqHA}
Let $A$ be a quasi-Hopf algebra such that its radical is a quasi-Hopf ideal,
and $A/\Rad A \cong \k[\Z_m]$ as algebras, for some $m\in \Bbb N$ not divisible by primes $\leq 7$.
Then $A$ is equivalent by a twist to one of the following:
\begin{enumerate}
  \item a radically graded finite-dimensional Hopf algebra $A$ such that
  $$ A/ \Rad A \cong \k [\Z_m] \mbox{, or}$$
  \item a semisimple quasi-Hopf algebra $\k[\Z_m]$, with associator given by $\omega_s \in H^3(\Z_m, \k^{\times})$, $s\in \{1,...,m-1\}$, or
  \item a quasi-Hopf algebra $A(H,s)$, where $H$ is a radically graded Hopf algebra such that $H/ \Rad H \cong \k [\Z_{m^2}]$, and $s \in \Upsilon(H)$.
\end{enumerate}
\end{thm}

We will give the proof in Subsection \ref{subsection:proofmainthm}.

Recall that a tensor category $\bC$ is \emph{pointed} if every simple object of $\bC$ is invertible. Invertible objects form a group. The
previous Theorem implies the corresponding statement for pointed finite tensor categories.

\begin{cor}
Let $\bC$ a pointed finite tensor category whose simple objects form a cyclic group of order $m$, where $m$ is not divisible by $2,3,5,7$. Then $\bC$ is equivalent to one of the following:
\begin{enumerate}
    \item the category of finite dimensional $H$-modules, for $H$ a radically graded finite-dimensional Hopf algebra such that $H/ \Rad H \cong \k [\Z_m]$, or
    \item a semisimple category $\Rep_{\omega_s}(\Z_m)$, or
    \item the category of finite dimensional $A(H,s)$-modules, for some radically graded Hopf algebra $H$ such that $H/ \Rad H \cong \k [\Z_{m^2}]$, and $s \in \Upsilon(H)$.
\end{enumerate}
\end{cor}
\bp
The fact that the category is pointed implies that its objects have integer Frobenius-Perron dimension, so by \cite{EO} it is the category of finite dimensional modules of some quasi-Hopf algebra $A$. This quasi-Hopf algebra is basic (because $\bC$ is pointed), so the result follows from Theorem \ref{thm:classificationqHA}.
\ep

The organization of this paper is the following. In Section \ref{section:preliminaries} we describe some tools which we use in the rest of the work. The two key results are the classification of pointed Hopf algebras over abelian groups given by Andruskiewitsch-Schneider, and the equivariantization procedure.

In Section \ref{section:gradedqHA} we construct basic radically graded quasi-Hopf algebras over $\Z_m$ as a generalization of the family $A(q)$ in \cite{G}. Using some methods in Etingof-Gelaki's works, we prove that these are all the basic radically graded quasi-Hopf algebras over $\Z_m$ up to twist equivalence.

After that, we consider liftings of these graded algebras in Section \ref{section:liftings}. We prove that each basic quasi-Hopf algebra whose associated radically graded quasi-Hopf algebra has trivial associator is a Hopf algebra, as in \cite{EG3}. For each non-semisimple basic radically graded quasi-Hopf algebra with non-trivial associator, we prove that any lifting $A$ can be extended to a Hopf algebra $H$ as in the graded case, so $\Rep H$ is the equivariantization of $\Rep A$ for some action of $\Z_m$; for an analogous procedure see \cite{EG4}. In this way we can describe such $A$ using the inverse procedure, the de-equivariantization of $\Rep H$ for an inclusion of $\Rep \Z_m$ (a result of Masuoka in \cite{Ma} reduces it to the graded case), and we complete the classification.

In Section \ref{section:classificationZpn}, we apply the previous classification to the case $m=p^n$ for some prime $p$ and some $n \in \N$. Such description is important for the general case, where we reduce some results to the case $p^n$.

\textbf{Acknowledgments.} The author's work was supported by CONICET,
and was done mainly during his visit to MIT in February-May 2009,
supported by Banco Santander Rio SA. The author thanks MIT for its warm hospitality, and especially his host
Professor Pavel Etingof for posing the problem, guidance and explanations about tensor categories,
and for many important suggestions that influenced this work.
He is also grateful to Cesar Galindo for many stimulating discussions.

\section{Preliminaries}\label{section:preliminaries}

For any Hopf algebra $H$, $\Delta$, $\epsilon$ and $S$ will denote the coproduct, counit and antipode, respectively. For the coproduct we will use Sweedler notation: for any $c \in C$, $\Delta(c)= c_1 \otimes c_2$.

For each tensor category $\bC$ we denote by $\bZ(\bC)$ the Drinfeld center of $\bC$. For each object $X$ of $\bC$, we denote by $FPdim \, X$ the Frobenius-Perron dimension of $X$, see \cite{EO}.

To begin with, we will describe some topics. First we give a brief introduction to Yetter-Drinfeld modules over a Hopf algebra $H$. Second, we consider the equivariantization and de-equivariantization procedures, for a better description see \cite{dgno},\cite{EG4} and \cite{ENO2}.

Also we consider the lifting theory for pointed Hopf algebras and the main results, see \cite{AS4}. For these Hopf algebras, we will give a brief characterization of their duals, which give place to basic Hopf algebras; i.e. their radical is a Hopf ideal, and $H/\Rad H \cong \Fun G$ for some finite group $G$.

\subsection{Yetter-Drinfeld modules and Drinfeld center of Hopf algebras.}

We recall the definition of a Yetter-Drinfeld module over a Hopf algebra in order to write the formulas defining this notion.

\begin{defn}
Let $H$ be a Hopf algebra. A left Yetter-Drinfeld module $M$ over $H$ is a left $H$-module $M$, with action denoted by $\cdot: H\otimes M \rightarrow M$, which is also a left $H$-comodule, with coaction $\delta:M \rightarrow H\otimes M$, $\delta(m)=m_{(-1)}\otimes m_{(0)}$, satisfying:
\begin{equation}\label{YDrelation}
    \delta(h\cdot m) = h_1m_{(-1)}S(h_3) \otimes h_2 \cdot m_{(0)}, \qquad h\in H, m \in M.
\end{equation}
Morphisms of Yetter-Drinfeld modules are $H$-linear morphisms which also preserve the comodule structure. The category of left $H$-comodules is denoted by $\hyd$: it is a tensor category, which inherits the action in the tensor product as $H$-modules, and coaction $\delta_{M\otimes N}: M\otimes N \rightarrow H\otimes M\otimes N$,
\begin{equation}\label{YDtensorcoaction}
    \delta_{M\otimes N}(m\otimes n)=m_{(-1)}n_{(-1)} \otimes m_{(0)} \otimes n_{(0)}, \qquad m\in M, \, n\in N.
\end{equation}
This category is braided, where the braiding for each pair $M,N \in \hyd$ is given by $c_{M\otimes N}: M\otimes N \rightarrow N\otimes M$,
\begin{equation}\label{YDbraiding}
    c_{M\otimes N}(m\otimes n)= m_{(-1)}\cdot n \otimes m_{(0)}, \qquad m \in M, n \in N.
\end{equation}
\end{defn}

\begin{rem}
\eqref{YDrelation} is equivalent to the following:
\begin{equation}\label{YDrelation2}
    (h_1\cdot m)_{(-1)}h_2 \otimes (h_1\cdot m)_{(0)} = h_1m_{(-1)} \otimes h_2 \cdot m_{(0)}, \qquad h\in H, m \in M.
\end{equation}
\end{rem}

The category $\hyd$ is equivalent to the category $\Rep(D(H))=\bZ(\Rep(H))$.

\subsection{Equivariantization and de-equivariantization.}

Let $\bC$ be a finite tensor category. Denote by $\aut\bC$ the category which have as objects the tensor auto-equivalences of $\bC$, and its morphisms are isomorphisms of tensor functors. It is a monoidal category, whose tensor product is the composition of tensor functors.

For any group $\Gamma$ denote by $\underline{\Gamma}$ the category whose objects are elements of $\Gamma$, its morphisms are just the identities on each object, and the tensor product corresponds to the multiplication in $\Gamma$.

\begin{defn}
An \emph{action} of a group $G$ on a finite tensor category $\bC$ is a monoidal functor $\bF:\underline{\Gamma} \rightarrow \aut \bC$.
\end{defn}

In this way, we have a collection of functors $\{F_g: g\in \Gamma\} \subset \aut \bC$, and isomorphisms
$$ \gamma_{g,h}: \xymatrix{ F_g \circ F_h \ar[r]^{\sim} & F_{gh} }, \quad g,h \in \Gamma, $$
defining the tensor structure of the functor $\bF$.

\begin{defn}
Let $\Gamma$ be a finite group acting on a finite tensor category $\bC$. A $\Gamma$\emph{-equivariant object} of $\bC$ is an object $X\in \bC$ with a family of isomorphisms $u_g: F_g(X) \rightarrow X$ such that for all pairs $g,h \in \Gamma$ the following diagram commutes:
$$ \xymatrix@C=.9in{ F_g(F_h( X )) \ar[r]^{F_g(u_h)}\ar[d]^{\gamma_{g,h}} & F_g(X) \ar[d]^{u_g} \\ F_{gh}(X) \ar[r]^{u_{gh}} & X.  } $$

A \emph{morphism of equivariant} objects $\beta: (X,(u_g)_{g\in \Gamma}) \rightarrow (Y,(v_g)_{g\in \Gamma})$ is a morphism $\beta:X \rightarrow Y$ in $\bC$ such that for all $g\in \Gamma$, $\beta \circ u_g=v_g \circ F_g(\beta)$. The category of $\Gamma$-equivariant objects is called the \emph{equivariantization} of $\bC$, and will be denoted $\bC^{\Gamma}$.
\end{defn}

For such category, we have a natural inclusion $\iota: \Rep \Gamma  \rightarrow \bC^{\Gamma}$.
\medskip

We consider the inverse procedure. Consider a finite tensor category $\bD$ such that $\bZ(\bD)$ contains a Tannakian subcategory $\Rep\Gamma$ for some finite group $\Gamma$, and the composition $\Rep \Gamma \rightarrow \bZ(\bD) \rightarrow \bD$ is an inclusion. The algebra $Fun(\Gamma)$ of functions $\Gamma \rightarrow \k$ is an algebra in the tensor category $\Rep \Gamma$: the group $\Gamma$ acts on $Fun(\Gamma)$ by left translations. In this way $Fun(\Gamma)$ is an algebra in the braided category $\bZ(\bD)$.

\begin{defn}
The category of $Fun(\Gamma)$-modules in $\bD$ is called the \emph{de-equivariantization} of $\bD$, and will be denoted by $\bD_{\Gamma}$. It is a tensor category.
\end{defn}

\medskip

We will use the following result about equivariantization and deequivariantization. For a complete reference and proofs about these facts, see \cite{dgno}.

\begin{thm}\label{thm:equiv-deequiv}
\vi Let $\Gamma$ be a finite group acting on a finite tensor category $\bC$. Then $\Rep \Gamma$ is a Tannakian subcategory of $\bZ(\bC^{\Gamma})$ (that is, the braiding of $\bZ(\bC^{\Gamma})$ restricts to the symmetric braiding of $\Rep \Gamma$), and the composition of $\Rep \Gamma \rightarrow \bZ(\bC^{\Gamma})$ with the forgetful functor $\bZ(\bC^{\Gamma}) \rightarrow \bC^{\Gamma}$ is the natural inclusion $\iota$.

\smallskip
\vii The procedures of equivariantization and deequivariantization are inverse to each other.
\end{thm}

\begin{exa}\label{example:semisimple-equiv}
We describe here an example over pointed semisimple categories. Although we shall work over non semisimple categories, we shall consider the semisimple part of some pointed ones and this will be useful in what follows.

Consider an action of a group $\Gamma$ over the category $\bC=Vec_{K,\omega}$, where $K$ is an abelian group and $\omega \in H^3(K, \k^{\times})$. We will denote the simple elements of $\bC$ just with the elements of $K$. We assume that the action over the objects is trivial; that is, $F_\gamma (X)=X$ for all object $X$ and all $\gamma \in \Gamma$. In this way, following the description on \cite[Section 7]{T} and using that the action is trivial on objects, the action is described by an element $\psi \in H^2(\Gamma,\hat{K})$:
$$ \psi(\gamma_1, \gamma_2): F_{\gamma_1} \circ F_{\gamma_2} \rightarrow F_{\gamma_1 \gamma_2}, \quad \gamma_i \in \Gamma.$$
From the tensor structure of each $F_{\gamma}$ we have an element $\xi \in  H^2(K,\hga)$,
$$ \xi(k_1,k_2)(\gamma): F_{\gamma}(k_1) \otimes F_{\gamma}(k_2)=k_1+k_2 \rightarrow F_{\gamma}(k_1+k_2)=k_1+k_2,$$
where $k_i \in K, \gamma \in \Gamma$.

We want to describe actions such that $\bC^{\Gamma}$ is pointed: it can be derived from \cite{N}, since by \cite{Ni} we have $\bC^{\Gamma} \cong (\bC \rtimes \Gamma)^*_{\bC}$. For our context, we derive that $\Gamma$ is abelian (notice also that we have an inclusion of $\Rep \Gamma$ in $\bC^{\Gamma}$) and so $\Rep \Gamma \cong Vec_{\hga}$. For a description as in \cite{N}, the action of $\hga$ on $K$ is trivial, and $\bC \rtimes \Gamma \cong Vec_{K \times \hga}$.

As $ FPdim (\bC^{\Gamma})= |\Gamma| FPdim \bC = |\Gamma| |K|$, $\bC^{\Gamma}$ has $|\Gamma| |K|$ non-isomorphic simple objects. Such objects are pairs $\left(k,(u_{\gamma})_{\gamma \in \Gamma}\right)$, for scalars $u_{\gamma} \in \k^{\times}$ satisfying $u_{\gamma_1}u_{\gamma_2}= \psi(\gamma_1,\gamma_2)(k) u_{\gamma_1 + \gamma_2}$. Therefore two simple objects $\left(k,(u_{\gamma})_{\gamma \in \Gamma} \right)$ and $\left(k,(v_{\gamma})_{\gamma \in \Gamma}\right)$ are related by an element $f \in \hga$ such that $v_{\gamma}=u_{\gamma}f(\gamma)$ for all $\gamma \in \Gamma$, which are isomorphic if and only if $f=1$. In this way we identify simple elements in $\bC^{\Gamma}$ as pairs $(k,f) \in K\times \hga$.

Also for any fixed $k$, there exist $|\Gamma|$ elements $\left(k,(u_{\gamma})_{\gamma \in \Gamma}\right)$, and
$$ \psi(\gamma_1, \gamma_2)(k)  = u_{\gamma_1} u_{\gamma_2}u_{\gamma_1+\gamma_2}^{-1}  = \psi(\gamma_2, \gamma_1)(k). $$
Therefore $\psi(\gamma_1, \gamma_2) = \psi(\gamma_2, \gamma_1)$ for all $\gamma_i \in \Gamma$, and from the relation given in \cite{T} we derive that $\xi(k_1,k_2)=\xi(k_2,k_1)$ for all $k_i \in K$. The elements of $H^2(K,\hga)$ parameterize central extensions of $K$ by $\hga$, and if $L$ is the corresponding to $\xi$, then $L$ is abelian and we can identify $\bC^{\Gamma}= Vec_{L,\widetilde{\omega}}$ for some $\widetilde{\omega} \in H^3(L,\k^{\times})$, because the tensor product in $\bC^{\Gamma}$ satisfies under the previous considerations:
$$ (k_1,f_1) \otimes (k_2 ,f_2)=(k_1+k_2,f_1+f_2+\xi(k_1,k_2) ), \qquad k_i \in K, f_i \in \hga. $$
Such $\widetilde{\omega}$ is the pullback of $\omega$ under the projection $\pi:L \rightarrow K$ corresponding to the extension, because the forgetful functor $\bC^{\Gamma} \rightarrow \bC$ is a tensor functor. This can also be derived from Naidu's work.

Also $\psi \in H^2(\Gamma, \hat{K})$ is the element corresponding to the dual extension $\hat{L}$ of $\Gamma$ by $\hat{K}$.
\smallskip

Note that, given a morphism $T: \hga \rightarrow \hat{L}$ such that for all $f_1,f_2 \in \hga$,
$$ \langle T(f_1), (0,f_2)\rangle = 1, $$
the function $\omega: K^3 \rightarrow \k^{\times}$ given by
$$ \omega(k_1,k_2,k_3)= \langle T\left( \xi(k_2,k_3) \right), (k_1,0) \rangle, \qquad k_i \in K, $$
defines an element in $H^3(K, \k^{\times})$, which will be denoted also by $\omega$. The pullback $\widetilde{\omega}$ of such element is trivial in $H^3(L,\k^{\times})$. Indeed, if $\alpha: L^2 \rightarrow \k^{\times}$ is the function
$$ \alpha \left( (k_1,f_1), (k_2,f_2) \right)=  \langle T(f_2), (k_1,f_1) \rangle = \langle T(f_2), (k_1,0) \rangle ,$$
then $\delta^2(\alpha)=\widetilde{\omega}$. In this way, $\bC^{\Gamma}\cong Vec_L$, and we have an inclusion $$ \Rep \Gamma \cong Vec_{\hga} \hookrightarrow Z(Vec_L) \cong Vec_{L \oplus \hat{L}}, $$
which composed with the forgetful functor to $\bC^{\Gamma}=Vec_L$ gives the canonical inclusion $\Rep \Gamma \hookrightarrow Vec_L$, so we have an inclusion of groups $\hga \hookrightarrow L \oplus \hat{L}$, which composed with the projection to the first component gives the inclusion $\hga \hookrightarrow L$.
\end{exa}

\begin{exa}\label{example:semisimple-deequiv}
Consider now the de-equivariantization of $\bD=Vec_L$, given by an inclusion of $\Rep \Gamma$ as a Tannakian subcategory of $Z(\bD)$, which factorizes through the center $Z(Vec_L) \cong Vec_{L \oplus \hat{L}}$; $\Gamma$ and $L$ are abelian groups as before, and we call $K$ the corresponding quotient group, which we also assume abelian. Therefore we have an inclusion $\iota: \hga \rightarrow  L$, and a morphism $T: \hga \rightarrow \hat{L}$, such that for all $f_1,f_2 \in \hga$, $ \langle T(f_1), (0,f_2)\rangle = 1$. Such $T$ parameterizes the natural morphisms $c_{V,-}: V\otimes - \rightarrow - \otimes V$ for each $V \in \Rep \Gamma$ viewed as an element of $\bD$.

Consider $L$ as an extension of $K$ by $\hga$, in such a way the inclusion $\iota$ is the canonical one, and it corresponds to an element $\xi \in H^2(K,\hga)$. The algebra $A=\Fun \Gamma$ is just the sum $A=\oplus_{f \in \hga} f $ as element of $Vec \hga$, with the canonical product, so we consider $A=\oplus_{f \in \hga} (0,f) $ inside $\bD$. By the previous considerations, we obtain $\bD_\Gamma \cong Vec_{K, \omega}$ for some $\omega \in H^3(K,\k^{\times}$.

The functor $F: \bD\rightarrow \bD_{\Gamma}$, $F(X)=A \otimes X$ is a monoidal functor, where the natural isomorphisms
$$ J_{X,Y}: F(X) \otimes_A F(Y) \rightarrow F(X\otimes Y) $$
are given by the natural isomorphisms induced by $T$ followed by the multiplication in $A$. Considering the monoidal functor axiom we deduce that $\omega(k_1,k_2,k_3)= \langle T\left( \xi(k_2,k_3) \right), (k_1,0) \rangle$.
\end{exa}

\subsection{Pointed Hopf algebras and liftings.}

We recall the Andruskiewitsch-Schneider Classification Theorem for pointed Hopf algebras over abelian groups whose order is divisible by primes greater than 7, and a result about their categories of comodules, due to Masuoka.

\begin{defn}[\cite{AS4}]
Let $\Gamma$ be an abelian group. A \emph{datum} of finite Cartan type over $\Gamma$,
$$ \bD=\bD \left(\Gamma, (g_i)_{i=1,...,\theta}, (\chi_i)_{i=1,...,\theta}, A=(a_{ij})_{i,j=1,...,\theta} \right), $$
consists of elements $g_i \in \Gamma$, $\chi_i \in \widehat{\Gamma}$ and a Cartan matrix of finite type $A$ satisfying for all $i,j$
$$ q_{ij}q_{ji}=q_{ii}^{a_{ij}}, \quad q_{ii}\neq 1, $$
where we define $q_{ij}:=\chi_j(g_i)$.
\end{defn}

Now call $\Phi$ the root system of the Cartan matrix $A$, $\bX$ the set of connected components of the corresponding Dynkin diagram and $\alpha_1,...,\alpha_{\theta}$ a set of simple roots; we write $i\sim j$ if $\alpha_i,\alpha_j$ are in the same connected component. For each $J\in \bX$, $\Phi_J$ denotes the root system of the component $J$.

Fix a datum $\bD$. For each $\alpha= \sum_{i=1}^{\theta} k_i\alpha_i \in \Phi^+$, we define
\begin{equation}\label{formulagchialpha}
    g_{\alpha}:=\prod_{i=1}^{\theta} g_i^{k_i}, \qquad \chi_{\alpha}:=\prod_{i=1}^{\theta} \chi_i^{k_i}.
\end{equation}

For our purposes, we consider $q_{ii}$ of odd order, and coprime with 3 if $\alpha_i$ belongs to a connected component of type $G_2$. In such case the order of $q_{ii}$ is constant on each connected component $J\in \bX$, and we define $N_J$ as the order of any $q_{ii}$.
\medskip

We introduce now two families of parameters. First we consider a family
$$\lambda= (\lambda_{ij})_{i,j \in \{1,..., \theta\}, i \nsim j}$$ of elements of $\k$ satisfying the condition:
\begin{equation}\label{lambdacondition}
    \mbox{if } g_ig_j=1  \mbox{ or }\chi_i\chi_j \neq \epsilon  ,\mbox{then } \lambda_{ij} =0.
\end{equation}
The second family is $\mu= (\mu_{\alpha})_{\alpha \in \Phi^+}$, which elements are also in $\k$, satisfying the condition:
\begin{equation}\label{mucondition}
    \mbox{if } g_{\alpha}^{N_J}=1  \mbox{ or }\chi_{\alpha}^{N_J} \neq \epsilon  ,\mbox{then } \mu_{\alpha} =0, \quad \forall\alpha \in \Phi_J^+,  J\in \bX.
\end{equation}

In \cite{AS4}, for any family $\mu$ and any $\alpha \in \Phi$, they introduce an element $u_{\alpha}(\mu) \in \k [\Gamma]$, which belongs to the augmentation ideal of $\k[g_i^{N_i}]$. An important fact for our work is that $u_{\alpha}(0)=0$ for all $\alpha \in \Phi^+$, where $\mu=0$ denotes the family which consists of all parameters equal to 0.

Also there exist elements $x_{\alpha}, \alpha \in \Phi^+$, which determine a PBW basis (see \cite{AS4} and the references therein).

\begin{defn}[\cite{AS4}]
The Hopf algebra $u(\bD,\lambda, \mu)$ is generated by $\Gamma$ and $x_1, \ldots, x_{\theta}$, with the following relations:
\begin{eqnarray}
  gx_ig^{-1} &=& \lambda_i(g)x_i, \quad i=1,...,\theta, \, g\in \Gamma; \label{groupaction}\\
  ad_c(x_i)^{1-a_{ij}}x_j &=& 0, \quad i\neq j, i \sim j; \label{qserre} \\
  ad_c(x_i)x_j &=& \lambda_{ij}(1-g_ig_j), \quad i<j, i \nsim j; \label{linkingrelation}\\
  x_{\alpha}^{N_J} &=& u_{\alpha}(\mu), \quad \alpha \in \Phi_J^+, \, J\in \bX. \label{powerrootlinfting}
\end{eqnarray}
\end{defn}

\begin{rem}
\begin{enumerate}
  \item In \cite{AS4} they prove that the algebra $u(\bD,\lambda, \mu)$ is a Hopf algebra, where the coproduct is defined by $\Delta(g)=g\otimes g$ for all $g\in \Gamma$, and $\Delta(x_i)=x_i \otimes 1+g_i \otimes x_i$. Its group-like elements are $ G\left(u(\bD,\lambda, \mu) \right) = \Gamma$.
  \item The graded case (trivial lifting) corresponds to $\mu=0, \lambda=0$.
\end{enumerate}
\end{rem}

\begin{thm}[\cite{AS4}]\label{thm:liftingsAS}
Let $H$ a finite dimensional pointed Hopf algebra, with group of group-like elements $\Gamma=G(H)$. Assume that the order of $\Gamma$ is not divisible by primes $\leq 7$. Then there exist a datum $\bD$ and families $\lambda, \mu$ such that $H\cong u(\bD,\lambda,\mu)$.
\end{thm}

\medskip

\begin{defn}[See \cite{AS2} and references therein]
Let $H$ be a bialgebra. A 2-cocycle on $H$ is a bilinear map $\sigma: H \times H \rightarrow \k$, which satisfies the following conditions
\begin{eqnarray}
\sigma(a_1,b_1) \sigma(a_2b_2, c) &=& \sigma(a, b_1c_1)\sigma(b_2,c_2), \\
  \sigma(a,1) &=& \sigma(1,a)=\epsilon(a),
\end{eqnarray}
for all $a,b,c \in H$.

Given an invertible (with respect to the convolution product) 2-cocycle, we define a new product on $H$ given by
$$ a \cdot_\sigma:= \sigma(a_1,b_1)a_2b_2 \sigma^{-1}(a_3,c_3), \qquad a,b \in H. $$
Then $H$ with this product, the same unit and the same coproduct structure is a new bialgebra. We denote it by $H_\sigma$. If $H$ is a Hopf algebra with antipode $S$, define
$$ S_ \sigma (a)= \sigma \left(a_1, S(a_2) \right) S(a_3) \sigma^{-1} \left( S(a_4), a_5 \right), \qquad a \in H.  $$
Then $S_\sigma$ is an antipode for $H_\sigma$, so $H_\sigma$ is a Hopf algebra.
\end{defn}

The following property of these liftings for coradically graded pointed Hopf algebras shall help us when we want to describe the category of representation of duals of pointed Hopf algebras.

\begin{thm}[\cite{Ma}]\label{thm:masuokacocycle}
Given a datum $\bD$ and families $\lambda, \mu$, the Hopf algebra $u(\bD,\lambda, \mu)$ is a cocycle deformation of the associated graded Hopf algebra $u(\bD,0, 0)$.
\end{thm}

\begin{rem}\label{rem:tensorequivwithgraded}
By the previous Theorem, the category of $u(\bD,\lambda, \mu)$-comodules is tensor equivalent to the category of $u(\bD,\lambda, \mu)$-comodules, see \cite{S1}.

Consider now a basic Hopf algebra $H$ such that such that $H/\Rad H \cong \Fun G$, where $G$ is an abelian group as in Andruskiewitsch-Schneider Classification Theorem. Denote by $H_0$ its associated radically graded Hopf algebra. Then $H^*$ is a pointed Hopf algebra isomorphic to some $u(\bD,\lambda, \mu)$, and its associated coradically graded Hopf algebra is $H_0^*$, which is isomorphic to $u(\bD,0, 0)$. Therefore, $\Rep H$ is tensor equivalent to $\Rep H_0$, because they are isomorphic to the categories of comodules over their corresponding duals.
\end{rem}

\subsection{Duals of pointed Hopf algebras}

Recall the following result:

\begin{prop}[\cite{B}]\label{prop:Beattidual} Let $\Gamma$ be a finite abelian group, and $V \in ^{\k \Gamma}_{\k \Gamma} \mathcal{YD}$, with basis $v_1,...,v_{\theta}$ where $v_i \in V^{\chi_i}_{g_i}$ for some $g_i \Gamma$, $\chi_i \in \widehat{\Gamma}$, such that $\bB(V)$ is finite dimensional. Then $H^* \cong \bB(W) \# \k\widehat{\Gamma}$, where we consider $W \in ^{\k \widehat{\Gamma}}_{\k \widehat{\Gamma}} \mathcal{YD}$ with a basis $w_i \in W^{g_i}_{\chi_i}$.
\end{prop}

\begin{rem}
The corresponding braiding matrices of $V$ and $W$ coincide: $(\chi_i(g_j))_{1 \leq i,j \leq n}$.
\end{rem}

We will describe duals of non-trivial liftings of Hopf algebras $u(\bD,\lambda, \mu)$ (see \cite{B} for case $A_1\otimes....\otimes A_1$, that is quantum linear spaces). Consider the coradically graded Hopf algebra $H_0=\bB(V) \# \k [\Gamma]$, for some abelian group $\Gamma$ and some $V\in \gyd$ such that $V$ is a diagonal braided vector space of Cartan type; if $A=(a_{ij})$ is the associated Cartan matrix of finite type, consider a basis $y_1,..., y_{\theta}$ with $y_i \in V^{\chi_i}_{g_i}$ ($g_i\in \Gamma, \chi_i \in \widehat{\Gamma}$) satisfying $\chi_i(g_j)\chi_j(g_i)=\chi_i(g_i)^{a_{ij}}$ for all $i \neq j$.

Liftings (in the coradical sense) $H$ are characterized as in Theorem \ref{thm:liftingsAS}, with the linking relations \eqref{linkingrelation} and the power root vector relations \eqref{powerrootlinfting}. Such $H$ has a basis $\{ hy: h \in \Gamma, y \in B\}$, where $B= \{ x_{\alpha_1}^{n_1}\cdots x_{\alpha_k}^{n_j}: \alpha_1 > \ldots >\alpha_k, \, 0 \leq n_j \leq N_{\alpha_j}-1 \}$ for a fixed order of $\Delta^+$, and $N_{\alpha}=N_J$ for each $\alpha \in J$, $J \in \bX$.

Define for each $\gamma \in \widehat{\Gamma}$ and each $y\in B$ the element $f_{\gamma,y} \in H^*$, which satisfies:
\begin{equation}\label{dualbasis}
f_{\gamma,y}(gy')=\gamma(g)\delta_{y,y'}, \qquad g\in \Gamma, y'\in B.
\end{equation}
In this way, $\{ f_{\gamma,y}: \gamma \in \widehat{\Gamma},y\in B \}$ is a basis of $H^*$. We call $x_i:=f_{\epsilon,y_i}$ for any $i\in \{1,...,\theta \}$, and identify $\gamma=f_{\gamma,1}$ for any $\gamma \in \widehat{\Gamma}$.

\begin{lem}\label{lemma:productdual}
Consider $H, x_1,...,x_{\theta}$ as above. Then,
\begin{enumerate}
  \item $\widehat{\Gamma} \cup \{x_1,...,x_{\theta} \}$ generate $H^*$ as algebra.
  \item $\Rad H^*$ is the ideal generated by $\{x_1,...,x_{\theta} \}$.
  \item in $H^*$, $\gamma x_i = \gamma(g_i)x_i \gamma$ for any $i\in \{1,...,\theta\}, \gamma \in \widehat{\Gamma}$.
\end{enumerate}
\end{lem}
\bp
(1) This follows from \cite[Lemma 2.1]{EG1}.
\medskip

(2) Remember that $Corad(H)^{\bot}=Rad(H^*)$, so the radical of $H^*$ is the ideal generated by $x_1,...,x_{\theta}$ and $H^*/Rad(H^*) \cong \widehat{\Gamma}$, by (1).
\medskip

(3) We calculate this explicitly for each $g\in \Gamma, z\in B$,
\begin{eqnarray*}
  (\gamma x_i)(gz) &=& (\gamma \otimes x_i)\Delta(gz)= (\gamma \otimes x_i)(gz \otimes g + gg_i\otimes z+ \widetilde{\Delta}(z)) \\
  &=& \gamma(gg_i) \delta_{z,y_i}= \gamma(g_i) (x_i \otimes \gamma)(gz \otimes g + gg_i\otimes z+ \widetilde{\Delta}(z)) \\
  &=& \gamma(g_i) (x_i \otimes \gamma)\Delta(gz) = \left( \gamma(g_i)x_i\gamma \right) (gz),
\end{eqnarray*}
where $\widetilde{\Delta}(z)$ is a sum of terms which first or second tensor term vanishes by applying $\gamma$.
\ep

\begin{lem}\label{lemma:coproductdual}
Consider $H, x_1,...,x_{\theta}$ as above. Call
\begin{equation}\label{grouplikedual}
    \Omega:= \left\{ \gamma \in \widehat{\Gamma}: \gamma\left(\lambda_i(g_i^{N_i}-1) \right)= \gamma\left(\alpha_{ij}(g_ig_j-1) \right)=0, \, 1 \leq i \nsim j \leq \theta  \right\}.
\end{equation}
\begin{enumerate}
  \item For the coproduct on $H^*$ we have:
\begin{eqnarray}\label{coproductdual}
    \Delta(x_i)-x_i\otimes 1-\chi_i \otimes x_i &\in& \Rad(H^*)\otimes \Rad(H^*),
    \\ \Delta(\chi)-\chi \otimes \chi &\in& \Rad(H^*)\otimes \Rad(H^*).
\end{eqnarray}
  \item $\Omega$ is the group of group-like elements of $H^*$.
\end{enumerate}
\end{lem}
\bp (1) Note that $\Delta(f)(x,y)=f(xy)$ for all $f\in H^*, x,y \in H$. Now for any $\chi \in \widehat{\Gamma}$ we note that
$$ \Delta(\chi) - \chi \otimes \chi \in \left( \widehat{\Gamma} \otimes \Rad{H^*} \right) \oplus \left( \Rad{H^*} \otimes \widehat{\Gamma} \right) \oplus \left( \Rad{H^*} \otimes \Rad{H^*} \right)  $$
(it is straightforward that $\Delta(\chi)$ contains $\chi \otimes \chi$ as the component in $\widehat{\Gamma} \otimes \widehat{\Gamma}$). Evaluating in $(g, g'y)$ and $(gy, g')$ for $g,g' \in \Gamma$, $y\in B$, we deduce that there are no components in $\widehat{\Gamma} \otimes \Rad{H^*} \oplus \Rad{H^*} \otimes \widehat{\Gamma}$, because $gyg'=q \, gg'y$ for some $q\in \k^{\times}$.

In a similar way we deduce the formula for $\Delta(x_i)$. Note that for each $z,z' \in B$, we express $zz'$ as a sum of elements of $B$ replacing $x_{\alpha}x_{\beta}$ by a sum of elements of $B$ of the same degree (where degree means length of words, i.e. viewing these in the tensor algebra of $V$), or where we replace some powers $x_{\alpha}^{N_{\alpha}}$ by $u_{\alpha}(\mu)$, or replace $x_ix_j$ by $\chi_j(g_i)x_jx_i+\lambda_{ij}(1-g_ig_j)$. So after reordering terms, $zz'$ is a sum of terms of the same degree in $B$, or of terms of less degree which have as factor $u_{\alpha}(\mu)$ or $\lambda_{ij}(1-g_ig_j)$. Then,
$$ \Delta(x_i)(gz, g'z')= x_i(q\,gg'zz')=0, \quad \forall (z,z') \neq (x_i,1), (1,x_i), $$
because $u_{\alpha}(\mu),\lambda_{ij}(1-g_ig_j) \in \ker \epsilon$. Also,
\begin{align*}
    \Delta(x_i)(gx_i, g')= x_i(\chi_i(g')gg'x_i)= \chi_i(g')   ,
\\  \Delta(x_i)(g, g'x_i)= x_i(gg'x_i)=1    .
\end{align*}
so we prove \eqref{coproductdual}.

\medskip
(2) It follows from the previous analysis about the expression of $zz',\, z,z' \in B$; see also \cite{B}.
\ep

\bigskip

\section{Basic Graded quasi-Hopf algebras}\label{section:gradedqHA}

In what follows, $m$ will denote a positive integer.

Consider a finite dimensional radically graded quasi-Hopf algebra: $H=\oplus_{i \geq 0} H[i]$, where $I:= \Rad H = \oplus_{i \geq 1} H[i]$, $I^k=  \oplus_{i \geq k} H[i]$. In such case, $H[0]$ is semisimple and $H$ is generated by $H[0]$ and $H[1]$ (Lemma 2.1, \cite[]{EG1}).

Observe that if $H$ is also basic, then $H[0]=\Fun(\Gamma)$ for some finite group $\Gamma$, where the associator (being in degree 0) corresponds to a class in $H^3(\Gamma, \k^{\times})$. Also, by \cite{S2}, $H[1]$ is a free module over $H[0]$. Consider now the case $\Gamma=\Z_{m}$.
\bigskip

Let $\sigma,\chi$ be generators of $\Z_{m}, \Z_{m^{2}}$, respectively, related by the condition $\chi^{m}=\sigma$ (considering the canonical inclusion $\Z_{m} \subseteq \Z_{m^{2}}$). Let $\{1_b: 0\leq b \leq m^{2}-1 \}$ the set of idempotents of $\k [\Z_{m^{2}}]$, defined by the condition $\chi 1_b=q^{b}1_b$ ($q$ a primitive root of unity of order $m^{2}$). Also, let $\{\uno_b: 0\leq s \leq m-1 \}$ the set of idempotents of $\k [\Z_{m}]$: as it is noted in \cite{G},
\begin{equation}\label{idemp}
    \sum_{0 \leq i \leq m-1} 1_{mi+s} = \uno_s, \qquad 0 \leq s \leq m-1.
\end{equation}
Also by \cite{G}, $\{\omega_s:  0\leq s \leq m-1\}=H^3(\Z_{m}, \C^{\times})$, where $\omega_s: (\Z_{m})^3 \rightarrow \k^{\times}$ is defined by
\begin{equation}\label{3-cocycle}
\omega_s(i,j,k) = q^{si(j+k-(j+k)')}, \qquad (i' \mbox{ denotes the remainder of the division by }m).
\end{equation}

In consequence, if $H$ is basic radically graded, the associator (being in degree zero) is
\begin{equation}\label{assoc}
    \Phi_s:= \sum_{i,j,k=0}^{m-1} \omega_s(i,j,k) \uno_i \otimes \uno_j \otimes \uno_k,
\end{equation}
for some $0\leq s \leq m-1$, which is trivial if and only if $s=0$.

Let $J_s= \sum_{i,j=0}^{m^{2}-1} c(i,j)^s 1_i\otimes1_j$, where $c(i,j):= q^{i(j-j')}$. As it is proved in \cite{G}, $J_s$ is invertible and satisfies:
\begin{equation}\label{conditionsJs}
(\varepsilon \otimes \id)(J_s)= (\id \otimes \varepsilon)(J_s)=1, \qquad \Phi_s=d\, J_s.
\end{equation}

\subsection{Quasi-Hopf algebras $A(H,s)$}

Given a radically graded Hopf algebra $H=\oplus_{n \geq 0} H(n)$ generated by a group like element $\chi$ of order $m^{2}$ and skew primitive elements $x_1,...,x_{\theta}$ satisfying \ref{skewprimitives}, $H=R \# \k [\Gamma]$, where $R \in \gyd$ the algebra of coinvariants.  If $\dim H$ is finite, $m^2$ does not divide $b_id_i$ (because $q^{b_id_i} \neq 1$).

We will define a quasi-Hopf algebra $A(H,s)$ for each $s \in \Upsilon(H)$ (recall the definition of $\Upsilon(H)$ given in Section \ref{section:introduction}), such that $A(H, s)/ \Rad A(H,s) \cong \k[\Z_{m}]$, with associator given by $\omega_s \in H^3(\Z_m,\k^{\times})$.

Consider the twist quasi-Hopf algebra $(H_{J_s}, \Delta_{J_s}, \varepsilon, \Phi_{J_s}, S_{J_s}, \alpha_{J_s}\beta_{J_s}, 1)$ and its subalgebra $A(H,s)$ generated by $\sigma:=\chi^{m}$ and $x_1,...,x_k$. Note that if $H$ is finite dimensional,
$$\dim A(H,s) = \dim H/ m= m \dim R.$$

\begin{prop}\label{prop:examples}
$(A(H,s), \Delta_{J_s}, \varepsilon, \Phi_{J_s}, S_{J_s}, \sigma^{-s}, 1)$ is a quasi Hopf algebra, which is not twist equivalent to a Hopf algebra.
\end{prop}
\bp
To simplify notation, we simply call $A=A(H,s)$. First of all, $\Phi_s \in A\otimes A\otimes A$. Using that $1_zx_i=x_i1_{z-d_i}$,
\begin{eqnarray*}
  \Delta_{J_s}(x_i) &=& \sum_{z,y=0}^{m^{2}-1} \frac{c(z,y)^s}{c(z-d_i,y)^s}q^{b_iy} x_i 1_{z-d_i} \otimes 1_{y} + \frac{c(z,y)^s}{c(z,y-d_i)^s}  1_{z} \otimes x_i 1_{y-d_i} \\
  &=&  \sum_{y=0}^{m-1} q^{b_iy} \left( \sum_{k=0}^{m-1} q^{mk(b_i-sd_i)} x_i \otimes 1_{y+km} \right)
  \\ &&+ \sum_{z,y=0}^{m^{2}-1} q^{z((y+d_i)'-d_i-y')} 1_{z} \otimes x_i 1_{y} \\
  &=&  \sum_{y=0}^{m-1} q^{b_iy} x_i \otimes \uno_{y} + \sum_{k=0}^{m-1} \left(\sum_{j=0}^{m-d_i'-1} q^{sk(d_i'-d_i)} \uno_{k} \otimes x_i \uno_{j} \right.
  \\ && \left. + \sum_{j=m-d_i'}^{m-1} q^{sk(d_i'+m-d_i)} \uno_{k} \otimes x_i \uno_{j} \right).
\end{eqnarray*}
Therefore $\Delta_{J_s}(x_i) \in A \otimes A$, and $(A, \Delta_{J_s}, \varepsilon_{J_s}, \Phi_{J_s})$ is a quasi bialgebra.

Now, $\alpha_{J_s}= \sum_{z=0}^{m^{2}-1} c(-z,z)^s 1_z, \, \beta_{J_s}= \sum_{z=0}^{m^{2}-1} c(z,-z)^s 1_z$, so
\begin{eqnarray*}
  \alpha_{J_s}\beta_{J_s} &=& \sum_{z=0}^{m^{2}-1} c(-z,z)^s c(z,-z)^s 1_z = \sum_{z=0}^{m^{2}-1} q^{smz} 1_z \\
  &=& \sum_{k=0}^{m-1} q^{smk} \uno_k = \left( \sum_{k=0}^{m-1} q^{k} \uno_k \right)^{ms} = \sigma^{-s}.
\end{eqnarray*}
Remember that $S(x_i)=-x_i \chi^{b_i}$, so
\begin{eqnarray*}
S_{J_s}(x_i) &=& \beta_{J_s} S(x_i) \beta_{J_s}^{-1}= -x_i \left( \sum_{y=0}^{m^{2}-1} \frac{c(y+d_i,-y-d_i)^s}{c(y,-y)^s}q^{-b_iy} 1_y \right) \\
  &=& -x_i \left( \sum_{y=0}^{m^{2}-1} q^{-b_iy+s(y+d_i)(m-(y+d_i)'+y+d_i)-sy(m-y'+y)} 1_y \right) \\
  &=& -x_i \left( \sum_{k,l=0}^{m-1} q^{-b_il+s(l+d_i)(m-(l+d_i)'+l+d_i)-slm+km(sd_i-b_i)} 1_{km+y} \right) \\
  &=& -x_i \left( \sum_{l=0}^{m-1} q^{-b_il+s(l+d_i)(m-(l+d_i)'+l+d_i)-slm} \uno_{y} \right) ,
\end{eqnarray*}
where we use again that $m$ divides $b_i-sd_i$.
\ep

\subsection{Radically graded quasi-Hopf algebras as subalgebras of twisted Hopf algebras.} We prove now that any radically graded quasi-Hopf algebra over $\Z_m$ looks like the quasi-Hopf algebras in the previous section. This fact gives us a characterization of all such quasi-Hopf algebras, in order to classify them.

\begin{thm}\label{thm:projection}
Let $A=\oplus_{n \geq 0} A[n]$ be a finite dimensional radically graded quasi-Hopf algebra over $\Z_m$, with associator $\Phi_s$ for some $s$ and $A[1] \neq 0$. Then, there exists a finite dimensional radically graded Hopf algebra $H$ as above, where $H= \bB(V) \# \Z_{m^2}$ for some Yetter-Drinfeld module $V$ over $\Z_{m^2}$, and a graded quasi-Hopf algebra epimorphism $\pi: A \twoheadrightarrow \bar{A}:= A(H,s)$, which is the identity restricted to degree 0 and 1.
\end{thm}
\bp This proof is similar to the one of Theorem 3.1 of \cite{EG1}. Decompose $A[1]= \oplus_{0 \leq r < m} H_r[1]$, where $$A_r[1]= \{x \in A[1]: \sigma x\sigma^{-1}=Q^rx\},$$
$Q=q^{m}$ a primitive root of unity of order $m$. Note that if $x \in A_r[1]$, $\uno_ix=x\uno_{i-r}$. Also, by \cite{EO}, we have that $H_0[1]=0$.

Let $\hat{A}$ be the tensor algebra of $A[1]$ over $A[0]$: it is a quasi-Hopf algebra, and we have a canonical surjective homomorphism $\pi_1: \hat{A} \twoheadrightarrow A$. Let $\gamma$ be the automorphism of $\hat{A}$ defined by
$$ \gamma|_{A[0]}= \id, \qquad \gamma|_{A_r[1]}= q^r\id. $$

Consider $L$ the sum of all quasi-Hopf ideals of $\hat{A}$ contained in $\oplus_{i \geq 2} \hat{A}[i]$. Therefore $\ker \pi \subseteq L$, and $\gamma(L)=L$, so $\gamma$ acts over $\bar{A}:=\hat{A}/L$. We define $\bar{H}$ as the quasi-Hopf algebra generated by $\bar{A}$ and a group-like element $\chi$, where $\chi^{m}=\sigma$ ($\chi$ has order $m^2$), and $\chi z\chi^{-1}=\gamma(z)$ for all $z\in H$. Note that $\Ad(\sigma)=\gamma^{m}$, so it is well defined, and $\chi$ generates a group isomorphic to $\Z_{m^{2}}$.

We consider the twist $H:= \bar{H}^{J^{-1}}$, which is a finite dimensional radically graded Hopf algebra. In such case, it is of the way $H=R\# \Z_{m^2}$, for some braided graded Hopf algebra $R$ in the category of Yetter-Drinfeld modules over $\Z_{m^2}$. We consider skew primitive elements $x_1,...,x_k \in H[1]$ which are eigenvectors of $\Ad(\chi)$:
$$    \chi x_i\chi^{-1}=q^{d_i}x_i, \quad   \Delta(x_i)=x_i \otimes 1 +\chi^{b_i}\otimes x_i , \qquad b_i,d_i \in \Z_{m^{2}}$$
Therefore, $\sigma x_i\sigma^{-1}=q^{d_jm}x_i$; as $H_0[1]=0$, $m \nmid d_j$.

If we denote $\bar{\Delta}$ the coproduct of $\bar{H}$, $\bar{\Delta}(x_i) \in \bar{A}\otimes \bar{A}$ because $\bar{A}$ is a quasi-Hopf subalgebra of $\hat{H}$. As $\bar{\Delta}(x_i)=J\Delta(x_i)J^{-1}$, we have
\begin{eqnarray*}
  \bar{\Delta}(x_i) &=& \sum_{z,y=0}^{m^{2}-1} \frac{c(z,y)^s}{c(z-d_i,y)^s}q^{b_iy} x_i 1_{z-d_i} \otimes 1_{y} + \frac{c(z,y)^s}{c(z,y-d_i)^s}  1_{z} \otimes x_i 1_{y-d_i} \\
  &=& \sum_{y=0}^{m^{2}-1} q^{sd_i(y'-y)+b_iy} x_i\otimes 1_{y} + \sum_{z,y=0}^{m^{2}-1} q^{sz(y'-d_i-(y-d_i)')} 1_{z} \otimes x_i 1_{y-d_i} \\
  &=&  \sum_{y=0}^{m-1} q^{b_iy} \left( \sum_{k=0}^{m-1} q^{mk(b_i-sd_i)} x_i \otimes 1_{y+km} \right)
  \\ && + \sum_{z,y=0}^{m-1} q^{sz((y+d_i)'-d_i-y)} \uno_{z} \otimes x_i \uno_{y}.
\end{eqnarray*}
The first summand belongs to $\bar{A}\otimes \bar{A}$, so $b_i \equiv sd_i (m)$.

Now, the braided graded Hopf algebra $R$ in the category of Yetter-Drinfeld modules over $\Z_{m^2}$ is generated in degree 1; call $V:=R[1]$. Therefore there exists an epimorphism of Hopf algebras $H \twoheadrightarrow \bB(V) \# \Z_{m^2}$, which induces by twisting and restriction (note that the kernel of such map is generated in degree $\geq 2$) a surjective morphism $\bar{A}=A(H, s) \twoheadrightarrow A(\bB(V)\#\Z_{m^2}, s)$. As both algebras have the same degree 0 and 1 parts and $\bar{A}$ has no proper quasi-Hopf ideals generated in degree $\geq 2$, such surjective map is an isomorphism, and $H=\bB(V) \# \Z_{m^2}$.
\ep
\bigskip

\subsection{Generation in degree 1}

In what follows, consider $m$ odd. Strictly speaking, we consider radically graded Hopf algebras, which are dual of coradically graded Hopf algebras.

Although $H$ and $H^*$ are of the same type (as groups, $\Z_m \cong \widehat{\Z_m}$ canonically), to be consistent with the notation we consider a braided vector space of diagonal type $W$ as above and fix a basis $x_1,...,x_{\theta}$, where $x_i \in W_{\chi_i}^{g_i}$ for some $g_i \in \Z_{m^2}$, $\chi_i \in \widehat{\Z}_{m^2}$, so the braiding matrix is $(\chi_i(g_j))_{1 \leq i,j \leq n}$. Call $\mathcal{X}$ the set of connected components of $A$. By Heckenberger's classification, on each connected component:
\begin{itemize}
  \item it is a braiding of Cartan type (see \cite{H}): there exists a Cartan matrix $A=(a_{ij})$ such that for all $i,j$, $q_{ii}^{a_{ij}}=q_{ij}q_{ji}$, or
  \item 3 divides $m$ and $V$ is of type \eqref{hatB2}, \eqref{hatB3} (see Section).
\end{itemize}
Such Nichols algebra is $\Z_{\theta}$-graded, where each $x_i$ has degree $e_i$.

\emph{Consider first the Cartan case.} Let $\Delta_+$ be the set of positive roots of $A$. We know that for each $\alpha \in \Delta_+$, there exists an element $x_{\alpha}\in \bB(V)$ of degree $\alpha$, such that the $x_{\alpha}$'s determine a PBW basis, with height $N_I$, determined by $I \in \mathcal{X}$ if $\alpha \in I$ (here is important that $2,3$ do not divide $n$); see \cite{AS3} and the reference therein. Moreover,

\begin{thm}[\cite{AS4}, Thm. 5.5, see also \cite{A}] The algebra $\bB(W)$ is presented by generators $x_1, \ldots, x_{\theta}$ and relations
\begin{eqnarray}
  ad_c(x_i)^{1-a_{ij}}x_j  &=& 0, \qquad \label{quantumserre} \\
  x_{\alpha}^{N_I} &=& 0.  \qquad \label{heightroot}
\end{eqnarray}
\end{thm}

Therefore, the algebra $\tilde{A}$ is generated by the same relations, and
\begin{equation}\label{groupelement}
    \sigma^m=0, \qquad \sigma x_i\sigma^{-1}=q^{d_im}x_i \, \, (i=1,..., \theta).
\end{equation}
if $\sigma$ denotes the generator of $\Z_m \cong \widehat{\Z}_m$ (because the multiplication is not changed by twisting).
\medskip

Our goal now is to prove that $\pi:A \twoheadrightarrow \bar{A}$ as above is really an isomorphism. In order to do that, we will prove that the relations \eqref{quantumserre} and \eqref{heightroot} hold in $A$ (relations \eqref{groupelement} are satisfied because $\pi$ is an isomorphism in degree 0 and 1, and a morphism of algebras).

\begin{prop}\label{prop:quantumserre}
Let $\pi:A\twoheadrightarrow \bar{A}$ be as in Theorem \ref{thm:projection}, with $A$ finite dimensional. Then, for all $i \neq j$, $ad_c(x_i)^{1-a_{ij}}x_j=0$ holds in $A$.
\end{prop}
\bp
Suppose that $z_{ij}:=ad_c(x_i)^{1-a_{ij}}x_j \neq 0$ in $A$. Then $x_i,x_j, z_{ij}$ are linearly independent (because they are linearly independent in $H$). By the previous construction, $\hat{A}= A(T(V)\#\Z_{m^2}, \omega_s)$, so we look a the coproduct in $A$ from the corresponding in $\hat{A}$ and projecting.

In $\hat{H}:=T(V)\#\Z_{m^2}$, $z_{ij}$ is skew primitive: $\Delta(z_{ij})=z_{ij} \otimes1+g^{(1-a_{ij})b_i+b_j} \otimes z_{ij}$, because $z_{ij}$ is primitive in $T(V)$. So the subalgebra $B$ generated by $a,x_i,x_j$ and $z_{ij}$ in $\hat{A}$ is a quasi-Hopf algebra, because $\Delta_{\hat{A}}(z_{ij}) =J\Delta(z_{ij})J^{-1}$. Applying Theorem \ref{thm:projection} for $B$, there exists a projection $B \twoheadrightarrow \bar{B}= A(\bB(V_1), s)$, so $\bar{B}$ is also finite dimensional. As the braiding of $V_1$ is independent of the basis for which is calculated (it is of diagonal type, see \cite{AS1}), we can calculate it with respect to the basis $y_1=x_i$, $y_2=x_j$, $y_3= z_{ij}$. Let $(Q_{st})_{s,t=1,2,3}$ the corresponding matrix. Using that $V$ is of Cartan type, we have
\begin{align*}
    &Q_{11}=q_{ii},  &Q_{12}Q_{21}=q_{ij}q_{ji},
    \\ &Q_{22}=q_{jj}   ,  &Q_{13}Q_{31}= q_{ii}^{2-a_{ij}}    ,
    \\ &Q_{33}=q_{ii}^{1-a_{ij}}, &Q_{23}Q_{32}=q_{ii}^{a_{ij}(1-a_{ij})}q_{jj}^2.
\end{align*}
This braiding is of Cartan type, or standard $\hat{B}_2 \times A_1$, or as in \eqref{hatB3}, because the order of the elements in the diagonal are odd and $\bar{B}$ is finite dimensional. In any case, there exists a matrix $(m_{st})$ as in \cite{A} associated with the braiding $(Q_{st})$. Therefore at least two vertices are not connected (there exist $s \neq t$ such that $m_{st}=m_{ts}=0$):
\begin{itemize}
  \item If $Q_{23}Q_{32}=1$, then $q_{jj}^2=q_{ii}^{a_{ij}(a_{ij}-1)}=q_{jj}^{a_{ji}(a_{ij}-1)}$, so $\ord q_{jj}$ divides $2-a_{ij}a_{ji}+a_{ji}$. This is a contradiction because $\ord q_{jj}$ is odd, greater than 1.
  \item If $Q_{12}Q_{21}=1$, then $a_{ij}=a_{ji}=0$ and $q_{ii}^{-m_{13}+2}=Q_{11}^{-m_{13}}Q_{13}Q_{31}=1$. The unique possibility is $m_{13}=3$, in which case $m_{23}=m_{32}=0$, but this contradicts the previous item.
  \item If $Q_{13}Q_{31}=1$, then $\ord q_{ii}$ divides $2-a_{ij}$, which cannot happen by a similar argument.
\end{itemize}
From this contradiction, $z_{ij}=0$ in $A$.
\ep

\begin{prop}\label{prop:heightroots}
Let $\pi:A\twoheadrightarrow \bar{A}$ be as in Theorem \ref{thm:projection}, with $A$ finite dimensional. Then, for all $\alpha \in I$, $I \in \mathcal(X)$, $x_{\alpha}^{N_I}=0$ holds in $A$.
\end{prop}
\bp
Following notation in \cite{AS3}, consider $\widetilde{\bB}(V)$ the algebra generated by $V$, where the $x_i$'s are primitive, and where relation \eqref{quantumserre} holds for all $i \neq j$: that is, consider the quotient of the tensor algebra $T(V)$ by the braided Hopf biideal generated by the quantum Serre relations. Call $H_1:= \widetilde{\bB}(V) \# \Z_{m^2}$. For $A_1:= A(H_1,s)$ we have a surjective map of algebras $A_1\twoheadrightarrow A$, because of Proposition \ref{prop:quantumserre}, which is of quasi-Hopf algebras because they have the same structure in degree 0,1 and they are generated by these components. So we have the following picture:
$$  \xymatrix{ H_1 \ar@{->>}[rd] & \ar@{->>}[l]\hat{H} \ar@{->>}[d] \\  & H } \quad \leftrightsquigarrow \quad \xymatrix{ \hat{A}=A(\hat{H}, \omega_s) \ar@{->>}[r] \ar@{->>}[d] \ar@{->>}[rd] & \ar@{->>}[ld]A_1=A(H_1, \omega_s) \ar@{->>}[d] \\ \bar{A}=A(\bar{H}, \omega_s) & \ar@{->>}[l] A. }  $$
Call $\mathcal{K}(V)$ the subalgebra generated by the $x_{\alpha}^{N_I}$ in $\widetilde{\bB}(V)$: by Proposition 4.7 in \cite{AS3}, it is a braided Hopf subalgebra of $\widetilde{\bB}(V)$. In consequence, by twisting and restriction, the algebra $\mathcal{K}$ generated by $a$ and $x_{\alpha}^{N_I}$ in $A_1$ is a quasi-Hopf subalgebra of $A_1$.

Suppose that at least one of the $x_{\alpha}^{N_I} \neq 0$ in $A$. Therefore, the subalgebra of $A$ generated by $a$ and $x_{\alpha}^{N_I}$, is a non zero quasi-Hopf subalgebra of $A$ (it is the image of $\mathcal{K}$). Consider then $x_{\alpha}^{N_I}$ a non zero element of minimal degree: it is a non zero primitive element because of the degree consideration. Therefore the subalgebra $\mathcal{K}'$ generated by $a$ and $x_{\alpha}^{N_I}$ in $A$ is finite dimensional, and admits a projection over a finite dimensional quasi-Hopf algebra $\mathcal{K}''= A(R\#\Z_{n^2}, \omega_s)$.

Looking at $x_{\alpha}^{N_I}$ as an element of $\widetilde{\bB}(V)$, we have $c (x_{\alpha}^{N_I} \otimes x_{\alpha}^{N_I}) = x_{\alpha}^{N_I} \otimes x_{\alpha}^{N_I}$, because $N_I$ is the order of $q_{\alpha}$, where $q_{\alpha}$ is the scalar such that $c (x_{\alpha} \otimes x_{\alpha}) = q_{\alpha} x_{\alpha} \otimes x_{\alpha}$ (it depends just on the $\Z_{\theta}$-graduation). Call $z= x_{\alpha}^{N_I}$. In $R$, as $\Delta$ is an algebra morphism (inside the category of Yetter-Drinfeld modules) and $\Delta(z)= z\otimes 1+1\otimes z$; inductively,
$$ \Delta(z^k)= \sum_{j=0}^k \binom{k}{j} z^j \otimes z^{k-j}  $$
(here we use that $c(z\otimes z)= z \otimes z$). As $R$ is finite dimensional (because $\mathcal{K}''$ is finite dimensional), there exists $k$ such that $z^k=0$. Considering the minimal one, we derive that $z=0$ because the field is of characteristic 0. From this, $R=0$, which contradicts that some $x_{\alpha}^{N_I}$ is non zero.
\ep

\emph{Consider now $V$ of standard type \eqref{hatB2}} (see \cite{A} for definition of standard braiding: we do not consider the case \eqref{hatB3} at the moment, because it does not appear for $\Z_{m}$ as we shall prove in Section \ref{nichols}). By \cite{A}, we know that $\bB(V)$ is presented by generators $x_1,x_2$ and relations:
\begin{eqnarray}
  x_1^3 &=& x_2^N=0, \label{powerB2} \\
  \left( ad_c(x_1)x_2 \right)^3 &=& \left( ad_c(x_1)^2x_2 \right)^{N'}=0 ,\label{powerrootsB2} \\
  ad_c(x_1)^3x_2 &=& ad_c(x_2)^2x_1 =0 , \label{qserreB2}\\
  \left[ad_c(x_1)^2x_2, ad_c(x_1)x_2 \right]_c &=& 0. \label{newqserreB2}
\end{eqnarray}
where $N, N'$ denote the order of $\zeta, \zeta\xi^{-1}$, respectively.

\begin{prop}\label{prop:standardrelations}
Let $A$ be a finite dimensional quasi-Hopf algebra such that as in Theorem \ref{thm:projection} we have $\pi: A \twoheadrightarrow \bar{A}= A(\bB(V) \# Z_{m^2},s)$ for such 2-dimensional standard braided vector space. Then \eqref{powerB2}, \eqref{powerrootsB2}, \eqref{qserreB2} and \eqref{newqserreB2} hold in $A$.
\end{prop}
\bp
The strategy to prove them is to consider algebras $H_1$ as in the proof of Proposition \ref{prop:heightroots} such that the corresponding $A_1= A(H_1, s)$ projects onto $A$, in order to obtain quasi-Hopf subalgebras of $A$, then apply Theorem \ref{thm:projection} and derive a contradiction if these relations are non zero.
\medskip

\vi Relations \eqref{powerB2} are easily proved, because $x_1^3, x_2^N$ are primitive in $T(V)$ as braided Hopf algebra in $\zmyd$.
\medskip

\vii The second relation on \eqref{qserreB2} holds as in Proposition \ref{prop:quantumserre}, because $q_{22}q_{21}q_{12}=1$ as in the Cartan case, and this is what is used in the proof of such Proposition. For the second, it is better to consider $H_1$ as the quotient by the Hopf ideal generated by $x_1^3$, because in such case $ad_c(x_1)^3x_2$ is skew-primitive by \cite[Lemma 5.7]{A}. In such case, we work as in Proposition \ref{prop:quantumserre}, considering the braiding matrix $(Q_{ij})_{i,j=1,2,3}$ with respect to $y_1=x_1, y_2=x_2, y_3= ad_c(x_1)^3x_2$:
\begin{align*}
    &Q_{11}=\xi,  &Q_{12}Q_{21}=\zeta^{-1},
    \\ &Q_{22}=\zeta   ,  &Q_{13}Q_{31}= \zeta^{-1}    ,
    \\ &Q_{33}=\zeta^{-2}, &Q_{23}Q_{32}=\zeta^{-1}.
\end{align*}
This diagonal braiding is not associated with a Nichols algebra of diagonal type, because all the vertices are connected but all the $Q_{ii}$'s have odd order. We have a contradiction, so \eqref{qserreB2} hold in $A$. In a similar way, left hand side of \eqref{newqserreB2}, which we call $y_3$, is skew-primitive by \cite[Lemma 5.9]{A}. Considering the braiding matrix $(Q_{ij})_{i,j=1,2,3}$ with respect to $y_1=x_1, y_2=x_2, y_3$:
\begin{align*}
    &Q_{11}=\xi,  &Q_{12}Q_{21}=\zeta^{-1},
    \\ &Q_{22}=\zeta   ,  &Q_{13}Q_{31}= \zeta^{-2},
    \\ &Q_{33}=\zeta^{-2}, &Q_{23}Q_{32}=\zeta,
\end{align*}
We have a contradiction again, so also \eqref{newqserreB2} holds in $A$.
\medskip

\viii Now consider $\widetilde{\bB}(V)$ the quotient of $T(V)$ by the braided Hopf biideal generated by all the relations except \eqref{powerrootsB2}, and $H_1 =\widetilde{\bB}(V) \# \Z_{m^2}$. As before, let $A_1= A(H_1, s)$. If $\bB(V)=T(V)/I(V)$, the ideal $I(V)$ is generated in consequence by \eqref{powerB2}-\eqref{newqserreB2}, so $\left( ad_c(x_1)x_2 \right)^3$ is primitive in $\widetilde{\bB}(V)$, because it belongs to the kernel of the surjection onto $\bB(V)$ and is of minimal degree. Therefore $\left( ad_c(x_1)x_2 \right)^3=0$ in $A$ by an analogous proof as in Proposition \ref{prop:heightroots}.

If now $\widetilde{\bB}(V)$ denotes the quotient of $T(V)$ by the braided Hopf biideal generated by all the relations except $\left( ad_c(x_1)^2x_2 \right)^{N'}=0$, in such algebra $\left( ad_c(x_1)^2x_2 \right)^{N'}$ is primitive, and again it implies that $\left( ad_c(x_1)^2x_2 \right)^{N'}=0$ in $A$.
\ep
\bigskip

\subsection{Classification.} With the previous results we can describe all the radically graded quasi-Hopf algebras over $\Z_m$ for $m$ odd. We summarize this in the following result.

\begin{thm}\label{thm:classificationgradedqHA}
Let $A$ be a radically graded finite dimensional quasi-Hopf algebra such that for some odd integer $m$, $$A/Rad(A)= \k[\Z_{m}].$$ Then $H$ is twist equivalent to one of the following quasi-Hopf algebras:
\begin{enumerate}
  \item radically graded Hopf algebras $A$ such that $A/\Rad(A) \cong \k [\Z_{m}]$,
  \item semisimple quasi-Hopf algebras $\k[\Z_m]$ with associator given by $\omega_s \in H^3(\Z_m,\k^{\times})$, for some $s\in \{0,1,...,m-1\}$,
  \item an algebra $A(H,s)$, for some radically graded Hopf algebra $H$ such that $A/\Rad(A) \cong \k [\Z_{m^2}]$, and some $s \in \Upsilon(H)$.
\end{enumerate}
\end{thm}
\bp
Given a radically graded quasi-Hopf algebra, if its associator is trivial, it corresponds to a Hopf algebra which dual is coradically graded with coradical $\Z_m$. But this family is self-dual.

Consider now radically graded quasi-Hopf algebras $A$ with non-trivial associator $\Phi_s$. If rank of $A[1]$ over $A[0]$ is zero, $A$ is semisimple, so $A= \k[\Z_m]$. If rank of $A[1]$ over $A[0]$ greater than 0, by Theorem \ref{thm:projection} there exists a coradically graded Hopf algebra $H$ with coradical $\Z_{m^2}$ and a projection of quasi-Hopf algebras $\pi: A \twoheadrightarrow A(H,s)$. If $H$ is of Cartan type, by Propositions \ref{prop:quantumserre} and \ref{prop:heightroots}, the relations defining $A(H,s)$ are satisfied in $A$, so $\pi$ is an isomorphism. If $H$ is not of Cartan type, by Heckenberger's classification of diagonal braidings \cite{H}, it is of standard type with some connected component of the generalized Dynkin diagram not of Cartan type, and by Proposition \ref{prop:standardrelations} relations defining $A(H,s)$ are satisfied in $A$, so again $\pi$ is an isomorphism. This completes the proof.
\ep

\bigskip

\section{Liftings of quasi-Hopf algebras over $\k[\Z_{m}]$}\label{section:liftings}

In this section, for any radically graded quasi-Hopf algebra $A_0$ with associator $\Phi_s$ such that $A_0/Rad(A_0) \cong \k[\Z_{m}]$ we consider the possible liftings: that is, all the non-semisimple quasi-Hopf algebras $A$ such that the associated graded quasi-Hopf algebra (with respect to the radical filtration) is $A_0$.

By the previous section, such $A_0$ are related with radically graded Hopf algebras $H_0$ such that $H_0/Rad(H_0) \cong \k[\Z_{m^{2}}]$: $A_0=A(H_0,s)$. We will relate the liftings $A$ of $A_0$ with liftings $H$ of $H_0$.

We will use the same denomination of deformation as in \cite{EG3}: a deformation of a map $f_0$ is a map $f$ obtained by adding terms in degree higher than some degree $d$.

\textbf{We restrict to the case $m$ not divisible by primes $\leq 7$:} at the moment pointed Hopf algebras over abelian groups are completely classified for those groups whose order is not divisible by $2,3,5,7$ (see \cite{AS4}).

\subsection{Lifting of quasi-Hopf algebras with trivial associator.}

We begin with quasi-Hopf algebras whose coradical is a quasi-Hopf ideal, such that the corresponding graded quasi-Hopf algebra is a Hopf algebra; i.e. the corresponding associator is trivial. Remember the following result:

\begin{prop}[\cite{EG3}]\label{prop:liftingtrivialassoc}
Let $A$ be a finite dimensional quasi-Hopf algebra whose radical is a quasi Hopf ideal and the corresponding graded algebra $A_0=gr(A)$ is a Hopf algebra. If $H^3(A_0^*, \k)=0$, then $A$ is twist equivalent to a Hopf algebra.
\end{prop}

Fix a radically graded Hopf algebra $A_0$  and consider a set of skew-primitive elements $x_i$ and a group-like element $\gamma$ as in Section \ref{section:gradedqHA}, which generate $A_0$ as an algebra.

Write $m=p_1^{\alpha_1} \cdots p_k^{\alpha_k}$, the decomposition of $m$ as product of primes, $p_i >7$ by hypothesis. Define
\begin{equation}\label{definitionVpi}
V(p_i):= \left\{ j\in \{1,...,\theta\}:  b_id_i \not\equiv 0 (p_i^{\alpha_i})  \right\}.
\end{equation}

As we consider finite dimensional Hopf algebras, $q^{b_id_i} \neq 1$, or equivalently $m$ does not divide $b_id_i$. Therefore, $\cup_l V(p_l)=\{ 1,...,\theta\}$. Also for each pair $i\sim j$, we have $a_{ij},a_{ji} \in \{1,2,3\}$ and $a_{ij}b_id_i \equiv a_{ji}b_jd_j (m)$, so $i \in V(p_l)$ iff $j\in V(p_l)$. In this way each $V(p_l)$ is a union of connected components of the Dynkin diagram associated with the diagonal braiding of $A_0$.

\begin{lem}
Let $p_l, V(p_l)$ as above. Then $|V(p_l)| \leq 2$.
\end{lem}
\bp
Let $A_{p_l}=(a_{ij})_{i,j \in V(p_l)}$ be the Cartan matrix obtained by restriction of $A$: it is another finite Cartan matrix. We have
$$ a_{ij}b_id_i \equiv b_id_j+b_jd_i \equiv a_{ji}b_jd_j (p_l^{\alpha_l}). $$
If $m'=\prod_{k\neq l}p_k^{\alpha_k}$, and $\bar{q}=q^{m'}$, it is a root of unity of order $p_l^{\alpha_l}$, then the braiding $(\bar{q}^{b_id_j})_{i,j \in V_l}$ is of finite Cartan type, associated to a braided vector space $W$ with a basis $\bar{x}_i$, such that we can fix a generator $\bar{\sigma}$ of $\Z_{p_l^{\alpha_l}}$ satisfying:
$$ \bar{\sigma}\bar{x}_i\bar{\sigma}^{-1}=\bar{q}^{d_i}\bar{x}_i, \quad \Delta(\bar{x}_i) = \bar{x}_i\otimes 1+ \bar{\sigma}^{b_i} \otimes \bar{x}_i. $$
Therefore it is one in Section \ref{nichols}, and $|V(p_l)| =\dim W \leq 2$
\ep

Now we state an analogous result to \cite[Thm. 1.3]{EG3}, and adapt the proof.

\begin{thm}\label{thm:trivialassoc}
Let $A$ be a finite dimensional quasi-Hopf algebra whose coradical is a quasi-Hopf ideal such that $A/\Rad(A) \cong \k[\Z_{m}]$, and the associated graded quasi-Hopf algebra $A_0=gr(A)$ is a Hopf algebra. Then $A$ is twist equivalent to a Hopf algebra.
\end{thm}
\bp
By Proposition \ref{prop:liftingtrivialassoc}, it is enough to prove that $H^3(A_0^*, \k)=0$. Note that $A_0^*$ is a coradically graded Hopf algebra with $G(A_0^*)=\Z_{m}$. By the results in Section \ref{classifnichols}, it is of the way $A_0^*= \Z_{m} \ltimes \bB(V))$, where $V$ is a braided vector space of Cartan type, and $q$ is a root of order $m^2$:
$$ H^\bullet(A_0^*, \k) = H^\bullet(\bB(V), \k)^{\Z_m}= H^\bullet(\mathfrak{u}_q^+,\k)^{\Z_m}. $$
The last equality is proved in \cite{EG3}, (although $\bB(V)$ and $\mathfrak{u}_q^+$ can be non-isomorphic as algebras, we have $H^\bullet(\bB(V), \k)= H^\bullet(\mathfrak{u}_q^+,\k)$).

In \cite{GK} they prove that $H^\bullet(\mathfrak{u}_q^+,\k)=\sum_{w\in W} \C \eta_w \otimes S(\mathfrak{n}_+)$, where $W$ is the Weyl group, each $\eta_w$ has degree the length of $w$ (that we will denote $\ell(w)$) and $S(\mathfrak{n}_+)$ is the symmetric algebra of the positive part of the associated Lie algebra sitting in degree 2.

Define $\rho:= \frac{1}{2} \sum_{\alpha \in \Delta_+} \alpha$. A generator $\sigma$ of $\Z_m$ acts trivially on $S(\mathfrak{n}_+)$, and by a scalar $\lambda_w$ on each $\eta_w$. Such scalar is
$$ \lambda_w:=q^{-\sum_i n_id_i} \qquad \mbox{where }\sum_i n_i \alpha_i = \gamma_w:=\sum_{\alpha \in \Delta_+: w(\alpha)< 0} \alpha=\rho-w(\rho).$$

As $\Z_m$ acts trivially in $S(\mathfrak{n}_+)$ and by $q^d_i \neq 1$ on $w=s_i$, in order to prove that $H^3(\mathfrak{u}_q,\k)^{\Z_m}=0$ it is enough to prove that $\lambda_w \neq 1$ for any $w$ such that $\ell(w)=3$. Write $w=s_{i_1}s_{i_2}s_{i_3}$, where at least two of them are different.
\smallskip

Assume first that $i_1,i_2,i_3 \in V(p_l)$ for some prime $p_l$ dividing $m$. Therefore two of them are equal, because such component belongs to one of the sets $V(p_l)$, and any of them has at most two elements by the above Lemma; assume $i_2=i_3$ so $V(p_l)=\{i_1,i_2\}$. Such set corresponds to a subdiagram of type $A_1\times A_1$, $A_2$, $B_2$ or $G_2$. The condition $\lambda_w=1$ is equivalent to $\sum_i n_id_i \equiv 0 (m)$, which we can consider just modulo $p_l^{\alpha_l}$. Using the characterization in Section \ref{nichols} and a computation analogous to the one in \cite[Prop. 5.1]{EG3}, we conclude that $\lambda_w \neq 1$ in this case.
\smallskip

Assume now that not all belong to the same $V(p_l)$. Then all the $i_j$ are different, and we fix by simplicity $i_j=j$. In this way there exists one of them which is in a different component: assume $3\in V(p_l)$ and $1,2 \notin V(p_l)$; i.e. $p_l^{\alpha_l}$ divides $b_1d_1$, $b_2d_2$, but it does not divide $b_3d_3$.

Write $b_i=p_l^{\beta_i}a_i$, $d_i=p_l^{\gamma_i}c_i$, where $p_l$ does not divide $a_i$, $c_i$. Then $\beta_1+\gamma_1$, $\beta_2+\gamma_2 \geq \alpha_l$, but $\beta_1+\gamma_1< \alpha_l$. Also, $b_1d_2+b_2d_1 \equiv 0 (p^l)$: it follows because $b_1d_2+b_2d_1 \equiv 0 (m)$ if $1 \nsim 2$, or because we consider $1,2 \notin V(p_l)$ if $1 \sim 2$. Therefore $\min\{\beta_1+\gamma_2, \beta_2+\gamma_1\} \geq \alpha_l$. From all this equations,
$$ \min \{\beta_1, \beta_2\}+\min \{\gamma_1, \gamma_2\} \geq \alpha_l. $$
Suppose now that $\lambda_w=1$. Therefore $d_3 \equiv -n_1d_1-n_2d_2 (m)$, where at least one of $n_1,n_2 \in \{1,2,3\}$ ($n_3=1$ because $3$ is in a different connected component of the Dynkin diagram). In this way, as $p_l^{\alpha_l}$ does not divide $d_3$ (because it does not divide $b_3d_3$) we deduce that $\min \{ \gamma_1,\gamma_2 \} \leq \gamma_3$. Also, as $3$ is not connected with $1,2$, we have $ b_1d_3+b_3d_1 \equiv b_2d_3+b_3d_2 \equiv 0 (p_l^{\alpha_l})$,
so $$ b_3d_3 \equiv -b_3(n_1d_1+n_2d_2) \equiv d_3 (n_1b_1+n_2d_2) (p_l^{\alpha_l}).$$
It follows that $\beta_3+\gamma_3 \geq \min \{\beta_1,\beta_2\}+\gamma_3 \geq \min \{\beta_1, \beta_2\}+\min \{\gamma_1, \gamma_2\} \geq \alpha_l$, which contradicts the hypothesis $3 \in V(p_l)$. Therefore $\lambda_w \neq 1$ also in this case.
\smallskip

From all these computations, $H^3(\mathfrak{u}_q,\k)^{\Z_m}=0$.
\ep

\subsection{Equivariantization of liftings of quasi-Hopf algebras.}
We want to obtain from each quasi-Hopf algebra $A$ which is a lifting of $A_0$, a Hopf algebra $H$ which is a lifting of $H_0$.

\begin{thm}\label{thm:equivariantization}
Let $A$ be a quasi Hopf algebra such that $Rad(A)$ is a quasi-Hopf ideal, and $gr(A)=A_0=A(H_0, s)$. There exists an action of $\Gamma=\Z_{m}$ on the category $\bC=\Rep (A)$ which fixes the simple elements of $\bC$, such that the equivariantization $\bC^{\Gamma}$ is tensor equivalent to $\Rep(H)$, for some Hopf algebra $H$. Such Hopf algebra is a lifting of $H_0$, and there exists an inclusion of quasi Hopf algebras $A \hookrightarrow H^J$, for some twist $J \in H\otimes H$.
\end{thm}
\bp
The first step is to construct $H$. The idea is to 'extend' $A$ as for the radically graded case following the steps in \cite{EG3}. We recall the main steps in order to see that such proof still holds in our context.

Consider the automorphism $\chi=S^2$ of the algebra $A$. Define $\Delta_{\chi}(x):= (\chi \otimes \chi)(\Delta(\chi^{-1}(x)))$ for each $x \in A$. By \cite{D}, Proposition 1.2, there exists a twist $K$ such that $\Delta_{\chi}(x)=K\Delta(x)K^{-1}$ for all $x\in A$. Call $K_0$ the degree zero part of $K$, which commutes with $\Delta_0(x)$ for all $x$ ($\Delta_0$ denotes the coproduct of $A_0$), so $K_0=1$, and hence $K=1+hdt$.

Note that \cite[Lemma 4.1]{EG3} holds in our setting: it uses just the fact that the character $\lambda$, which determines the isomorphism of tensor functors $V \rightarrow \lambda \otimes V^{****} \otimes \lambda^{-1}$ in $\Rep A$ has order $m$. Therefore we conclude that $S^{2m}$ is an inner automorphism: there exists $b=\sigma+hdt \in A$ ($\sigma \in A_0$ is the fixed group element) such that $S^{2m}(x)=bxb^{-1}$ for all $x\in A$. Also \cite[Lemma 4.2]{EG3} applies here (it uses the fact that $A$ is a lifting of $A_0$ with $A_0/Rad(A_0)=\Z_{m}$ but not the particular structure of $A_0$), and then we can choose $b$ satisfying the relation
\begin{equation}\label{conditionb}
   K^{-1} \chi^{\otimes 2}(K^{-1})\cdots (\chi^{n-1})^{\otimes 2} (K^{-1}) = \Delta(b)(b^{-1} \otimes b^{-1}).
\end{equation}
Using the construction for the semidirect product explained in \cite[Section 3]{EG3}, we can define the quasi-Hopf algebra
$$ \wH:= (\k[\chi,\chi^{-1}] \ltimes A)/(\chi^n -b). $$
Note that this algebra is characterized by the multiplication of $A$, $\chi a\chi^{-1}= \chi(a)=S^2(a)$ and the relation $\chi^n=b$ established by the quotient.

Considering a lifting $J \in \wH \otimes \wH$ of $J_s$, $H:= \wH^{J^{-1}}$ is a quasi-Hopf lifting of $H_0$. As it is showed \cite[Theorem 1.3]{EG3} (we use really a generalization of this proof for $\Z_{m}$ in place of $\Z_p$, see proof of Theorem \ref{thm:trivialassoc}), $H^3(H^*_0, \k)=0$ when $H_0$ corresponds to Nichols algebras of Cartan type for $\Z_{m}$, so as in \cite[Theorem 4.3]{EG3} we can change $J$ for $JF$ for some $F=1+hdt$ such that $H$ is a Hopf algebra (still a lifting of $H_0$).
\medskip

Note that $\Rep \wH \cong \Rep H$. To complete the proof, define the action of $\Z_m$ on $\bC=\Rep A$ (this is analogous to the proof of \cite[Thm. 4.2]{EG4}): we call $h$ a generator of $\Z_m$, to distinguish it from $\chi \in \wH$. To do this, we have to define a collection of functors $\{F_k:=F_{h^k} \}_{k=0,1...,m-1} \subset \Aut(\bC)$. For each $(V,\pi_V)\in \Rep A$, consider $F_k(V)=V$, and $\pi_{F_k(V)}(a)=\pi_V(S^2(a))$ for all $a\in A$.

The natural isomorphism $\gamma_{i,j}: F_i(F_j(V)) \rightarrow F_{(i+j)'}(V)$ is given by the action $b^{\frac{(i+j)'-i-j}{n}} \in A$: explicitly, $F_k \circ F_j=F_{k+j}$ if $j+k <n$, and $F_k\circ F_j $ is related with $F_{j+k-n}$ up to the action of $b^{-1}$.

For this action, a $\Gamma$-equivariant object of $\bC$ is an object $X\in \bC$ together with a collection of linear isomorphisms $u_k: F_k(X)=X \rightarrow X$ such that
\begin{align*}
    u_k(a\cdot v) &= S^{2k}(a)u_k(v), \qquad a\in A, v \in V; \\
    u_k u_j &=u_{(k+j)'}b^{\frac{(i+j)'-i-j}{n}}.
\end{align*}
These relations are exactly the ones defining $H$ as we have seen if $u_k=\chi^k$, so we have an equivalence of categories between $\bC^{\Gamma}$ and $\Rep(H)$; the tensor product of these representations is the same as for representations of $H^{J}=\wH$. So this completes the proof.
\ep

\subsection{De-equivariantization of $\Rep(H)$ for liftings of $H_0$.}

We want to obtain $A$ from $H$ for each lifting $H$ of a Hopf algebra $H_0$ as above;this is possible thanks that the de-equivariantization procedure is the inverse of equivariantization. But then we want to know for which liftings $H$ of $H_0$ we can apply it, in order to obtain all the quasi-Hopf liftings $A$. That is, we want to know all inclusions $\Rep(\Z_m) \rightarrow \mathcal{Z}(\Rep(H)) \cong \Rep(D(H))$ such that they factorize the inclusion $\Rep(\Z_m) \rightarrow \Rep(H)$. We begin characterizing such functors.

Consider a radically graded Hopf algebra over $\Z_{m^2}$ such that $\Upsilon(H_0) \neq \emptyset$. We prove now a technical lemma which we need in what follows.

\begin{lem}\label{lemma:grouplikeorderm}
Fix $H$ a lifting of $H_0$ as in Section \ref{section:gradedqHA}. Then, $\sigma=\chi^m \in G(H)$
\end{lem}
\bp
To prove that $\chi^m$ is a group-like element is equivalent to prove that $\chi^m(g_ig_j)=\chi^m(g_{\alpha}^{N_{\alpha}})=1$ for each pair $i,j$ such that $\lambda_{ij}=0$ and for each positive root $\alpha$ such that $\mu_{\alpha} \neq 0$, by Lemma \ref{lemma:coproductdual}.

Consider $i \nsim j$ such that $\lambda_{ij} \neq 0$. Therefore $\epsilon= \chi_i\chi_j= \chi^{d_i+d_j}$, so $d_i+d_j \equiv 0 (m^2)$. Now for $s \in \Upsilon(H_0)$, $b_k \equiv sd_k (m)$ for all $k$. Then $b_i+b_j \equiv s(d_i+d_j) \equiv 0 (m)$, so
$$ \chi^m(g_ig_j)=\chi^m(g^{b_i+b_j})=q^{m(b_i+b_j)}=1. $$

Now consider a positive root $\alpha=\sum n_j\alpha_j$ such that $\mu_{\alpha}\neq 0$. Therefore
$$ \epsilon= \chi_{\alpha}^{N_{\alpha}}= \left( \prod_{j=1}^{\theta}\chi_j^{n_j} \right)^{N_{\alpha}} =\chi^{N_{\alpha}(\sum n_id_i)}. $$
Then $N_{\alpha}(\sum n_id_i) \equiv 0 (m^2)$. Consider $s$ as above, so we have
$$ N_{\alpha}(\sum n_ib_i) \equiv N_{\alpha}s(\sum n_id_i) \equiv 0 (m), $$
which implies that $\chi^m(g_{\alpha}^{N_{\alpha}})=1$.
\ep

\begin{prop}\label{prop:inclusionZm}
Fix $H$ a lifting of $H_0$ as in Section \ref{section:gradedqHA}. There is a bijection between:
\begin{enumerate}
  \item functors $F:\Rep(\Z_m) \rightarrow \mathcal{Z}(\Rep(H))$ such that $\Rep(\Z_m)$ is a Tannakian subcategory of $\mathcal{Z}(\Rep(H))$, and the composition $\Rep(\Z_m) \rightarrow \mathcal{Z}(\Rep(H)) \rightarrow \Rep(H)$,
  \item integers $s \in \Upsilon(H_0)$.
\end{enumerate}
\end{prop}
\bp
As before, $\chi$ denotes a generator of $\widehat{\Z_{m^2}}=\Z_{m^2}$, which satisfies $\chi(g)=q$ for our fixed root of unity of order $m^2$.

Consider a functor as in (1): it is given by the projection $H \twoheadrightarrow H/\Rad(H) \cong \k[\Z_{m^2}] \twoheadrightarrow \k[Z_m]$, where both projections are the canonical ones. In this way, we have the element $\gamma \in H \setminus \Rad(H)$, which is a preimage of the generator of $\widehat{\Z_{m^2}}=\Z_{m^2}$, as we defined in Section \ref{section:preliminaries}.

$\Rep \Z_m$ is semisimple, and we essentially have to identify the simple $\Z_m$-modules $M_i= \k v_i$ (we fix a non-zero vector $v_i$ of this one-dimensional vector space), $i=0,1,...,m-1$. By the equivalence $\bZ(\Rep H) \cong \hyd$, we consider $M_i \in \hyd$, where the action should be given by $\gamma \cdot v_i=q^{mi}v_i$.

We have to define an structure of $H$-comodule for each $M_i$, $\delta:M_i \rightarrow H\otimes M_i$. As $\dim M_i=1$, it is determined by a group-like element $\chi_i \in H$ for each $i=0,1,...m-1$, such that $\delta(v_i)=\chi_i \otimes v_i$.

In $\Rep G$, $M_i \otimes M_j \cong M_{(i+j)'}$, and we want a tensor inclusion. By \eqref{YDtensorcoaction}, this means that $\chi_i\chi_j=\chi_{(i+j)'}$, so $\chi_1=\chi^{ms}$ for some $s\in \{0,1,...,m-1\}$, and this determines $\chi_i=\chi^{msi}$ for all $i=0,1,...,m-1$. By Lemma \ref{lemma:grouplikeorderm}, all the $\chi^{msi}$ are group-like elements.

These action and coaction should satisfy \eqref{YDrelation2}. As $\chi, x_1,...,x_{\theta}$ generates $H$ as algebra, it is enough to prove this relation for these generators. We use here Lemma \ref{lemma:coproductdual}. When $h=\chi$ and $m=v_i$, as $\Rad H$ acts by 0, both sides of \eqref{YDrelation2} are equal to $q^{mi}\, \chi^{msi+1} \otimes v_i$. When $h=x_k$ and $m=v_i$, as $x_k$ acts by 0, the left and the right-hand sides of \eqref{YDrelation2} are, respectively,
\begin{align*}
    & (x_k\cdot v_i)_{(-1)}1 \otimes (x_k \cdot v_i)_{(0)}+ (\chi^{b_k}\cdot v_i)_{(-1)}x_k \otimes (\chi^{b_k} \cdot v_i)_{(0)}= q^{b_k mi} \, \chi^{msi}x_k \otimes v_i,
\\ & x_k \chi^{msi} \otimes 1 \cdot v_i + \chi^{b_k}\chi^{msi} \otimes x_k \cdot v_i = q^{d_k s mi} \, \chi^{msi}x_k \otimes v_i.
\end{align*}
Therefore the action satisfies \eqref{YDrelation2} if and only if $b_k \equiv sd_k (m)$ for all $k$.

Also the braiding of $\hyd$ restricts to the canonical symmetric braiding of $\Rep \Z_{m}$. In fact, for each pair $k,j$, the braiding $c_{M_k,M_j}: M_k \otimes M_j \rightarrow M_j \otimes M_k$ is, by \eqref{YDbraiding},
$$ c(v_k \otimes v_j)= (v_k)_{(-1)} \cdot v_j \otimes (v_k)_{(0)}= \chi^{msk} \cdot v_j \otimes v_k = q^{m^2skj} v_j \otimes v_k=v_j \otimes v_k. $$

\medskip
Reciprocally, consider $s \in \Upsilon(H_0)$. Define $\bF_s:\Rep \Z_m \rightarrow \hyd$ as the functor which is the induced by the projection $H \twoheadrightarrow H/\Rad(H) \cong \k[\Z_{m^2}] \twoheadrightarrow \k[Z_m]$ as modules, and for each $M_i= \k v_i$ define as before $\delta: M_i \rightarrow H\otimes M_i$, $\delta(v_i)=\chi^{msi} \otimes x_i$. By the previous computations, these structures satisfy the compatibility condition \eqref{YDrelation2}, so $\bF_s$ sends objects to objects. For morphisms, the semisimplicity of $\Rep \Z_m$ gives a canonical definition of $\bF_s$, preserving the abelian structures of categories.

As above $\bF_s$ is tensorial, and moreover is braided, if we consider the canonical symmetric braiding of $\Rep \Z_{m}$. So the proof is completed.
\ep

\begin{lem}\label{lemma:sneqt}
Suppose that $s \in \Upsilon(H_0)$. Then the de-equivariantization $(\Rep H)_{\Z_m}$ induced by the inclusion $\bF_s$ is $\Rep A$ for some basic quasi-Hopf algebra such that $A/\Rad A \cong \k[\Z_m]$ with associator given by $\omega_s \in H^3(\Z_m,\k^{\times})$.
\end{lem}
\bp
By \cite[Corollary 4.27]{dgno}, the category $(\Rep H)_{\Z_m}$ is integral, so it corresponds to $\Rep A$ for some quasi-Hopf algebra $A$. We apply the computations in Example \ref{example:semisimple-deequiv} to the semisimple part of these categories, and we obtain that the semisimple part of $\Rep A$ is $Vec_{\Z_m, \omega}$, where $\omega$ is given as in such example. As here we have $T: \Z_m \rightarrow \widehat{\Z}_{m^2}$, $T(a)(b)= q^{sab}$ for $q$ a root of unity of order $m^2$ as above, and the 2-cocycle defining $\Z_{m^2}$ as an extension of $\Z_m$ by $\widehat{\Z}_m \cong \Z_m$ is $$ \xi: \Z_m \times \Z_m \rightarrow \Z_m \cong < m > \subseteq \Z_{m^2}, \quad \xi(j,k)= (j+k)'-j-k, $$
we deduce that $\omega=\omega_s$.
\ep

In order to classify all the liftings of quasi-Hopf algebras, we have to classify all the possible inclusions $\Rep \Z_m \rightarrow \Rep D(H)$ for liftings $H$ of Hopf algebras which satisfy conditions in Section \ref{section:gradedqHA}, and consider their de-equivariantizations.

\begin{lem}\label{lemma:de-equivgradedcase}
Let $H_0$ be a radically graded Hopf algebra such that $H_0/\Rad H_0 \cong Z_{m^2}$. For each integer $s \in \Upsilon(H_0)$, the de-equivariantization of $\Rep H_0$ corresponding to the functor $\bF_s$ is $\Rep A(H_0,s)$.
\end{lem}
\bp
By the proof of Theorem \ref{thm:projection} and Theorem \ref{thm:classificationgradedqHA}, we can extend each $A(H_0,s)$ to $H^{J_s}$ in such a way we obtain $\Rep H_0$ as a equivariantization of $\Rep A(H_0,s)$ by an action of $\Z_m$ fixing the invertible elements, see the proof of Theorem \ref{thm:equivariantization}. By Theorem \ref{thm:equiv-deequiv}, each $\Rep A(H_0,s)$ is in consequence a de-equivariantization of $\Rep H_0$ by an inclusion of $\Rep \Z_m$. The result follows by the previous Lemma.
\ep

\subsection{Proof of Theorem \ref{thm:classificationqHA}.}\label{subsection:proofmainthm}

For $A$ a quasi-Hopf algebra as in the Theorem, consider its associated radically graded quasi-Hopf algebra $A_0$. If $A_0$ has trivial associator, then it is a Hopf algebra, and $A$ is twist equivalent to a Hopf algebra by Theorem \ref{thm:trivialassoc}. The dual algebra has coradical isomorphic to $\Z_m$, because dualizing the radical filtration we obtain the coradical filtration, see \cite{Mo}. By Theorem \ref{thm:liftingsAS}, its dual is a lifting $u(\bD, \lambda,\mu)$ for a datum $\bD$ over $\Z_m$, but by Theorem \ref{thm:masuokacocycle} it is a cocycle deformation of $A_0^*= u(\bD,0,0)$, so $A$ is twist equivalent to $A_0$, which is also of type $u(\bD',0,0)$ for some datum $\bD'$ over $\Z_m$.
\smallskip

Consider now the case when $A_0$ has non-trivial associator: by Theorem \ref{thm:classificationgradedqHA}, it is semisimple with non-trivial associator, or it is of the way $A(H_0,\omega_s)$ for some radically graded (and in consequence also coradically graded) Hopf algebra $H_0$ with group of group-like elements $\Z_{m^2}$ and $s \in \Upsilon(H_0)$. In the first case we are done, so consider the second. By Theorem \ref{thm:equivariantization}, the category $\Rep A$ admits an action of $\Z_m$ whose equivariantization is $\Rep H$, where $H$ is a lifting (in the radical sense) of the Hopf algebra $H_0$.

On the other hand, $\Rep H$ is tensor equivalent to $\Rep H_0$ by Remark \ref{rem:tensorequivwithgraded}, and this tensor equivalence induces an equivalence between the corresponding centers, which commutes with the forgetful functors. So an inclusion of $\Rep \Z_m$ in $\Rep H$ factorizing through the center of $\Rep H$ corresponds univocally to an inclusion in $\Rep H_0$ factorizing through the center of $\Rep H_0$, and that tensor equivalence induces also a tensor equivalence between the $A=\Fun \Z_m$-modules on such categories. That is, de-equivariantizations of $\Rep H$ are tensor equivalent to de-equivariantizations of $\Rep H_0$.

In consequence we reduce the problem to the graded case, and $\Rep H_0$ admits as many inclusions $\Rep G \rightarrow \Rep D(H)$ as numbers $s$ are in $\Upsilon(H_0)$. But each $s$ corresponds to a de-equivariantization $Rep A(H_0,\omega_s)= (\Rep H)_{\Z_m}$ by Lemma \ref{lemma:de-equivgradedcase}, so they correspond to all the de-equivariantizations. As these procedures are inverse one of the other, $\Rep A \cong \Rep A(H_0,\omega_s)$ for some $s$, and in consequence $A$ is equivalent to $A(H_0,\omega_s)$.

\bigskip

\section{Explicit description of quasi-Hopf algebras over $\Z_{p^n}$, $p$ prime}\label{section:classificationZpn}

As an example of the previous result, we will describe all the basic finite-dimensional quasi-Hopf algebras $A$ such that $A/\Rad A \cong \k [\Z_{p^n}]$, for $p$ a prime greater than 7 and any $n\in\N$. This is based in the classification of pointed Hopf algebras over $\Z_{p^n}$, so first of all we describe all the possible Nichols algebras of finite dimension over $\k[\Z_{p^n}]$. It is done as in \cite{AS1} for $n=1$, and we shall obtain here an analogous description for the general case. Moreover, we can classify radically graded quasi-Hopf algebras over $\Z_{p^n}$ for any odd prime $p$ (by the general Theorem \ref{thm:classificationgradedqHA}). Then we restrict our attention to the case $p >7$, because of Theorem \ref{thm:liftingsAS}.

\subsection{Nichols algebras over $\Z_{p^n}$}\label{nichols}

We consider which are the possible Nichols algebras over $\Z_{p^n}$ following the description in \cite{AS1}. We fix $q$ a primitive root of unity of order $p^n$ and $g$ a generator of $\Z_{p^n}$, $p$ odd.

As in such work, we consider a basis $x_1,...,x_l$, where $x_i \in V_{g_i}^{\chi_i}$ for some $g_i \in \Z_{p^n}$ and some characters $\chi_i$. The characters are determined by $\chi_i(a)=q^{d_i}$, and we write $g_i=g^{b_i}$. Consider $a_i,c_i \in \N$ not divisible by $p$, and $\alpha_i, \gamma_i \geq 0$ such that $b_i=p^{\alpha_i}a_i$, $d_i=p^{\gamma_i}c_i$.

By \cite{H}, the braiding matrix $(q_{ij}:=\chi_i(g_j))$ is of Cartan type, or $p=3$ and the braiding is one of the following:
\begin{align}\label{hatB2}
\hat{B}_2: &\xymatrix{\circ^{\xi} \ar@{-}[r]^{\zeta^{-1}} & \circ^{\zeta}}, \qquad \xi \in \mathbb{G}_3, \zeta \in \k^{\times} \setminus \{ 1, \xi, \xi^2 \},
\\& \xymatrix{ \circ^{\zeta} \ar@{-}[r]^{\zeta^{-1}} & \circ^{\zeta} \ar@{-}[r]^{\zeta^{-1}} & \circ^{\zeta^{-3}} }, \quad \xymatrix{ \circ^{\zeta} \ar@{-}[r]^{\zeta^{-1}} & \circ^{\zeta^{-4}} \ar@{-}[r]^{\zeta^4} & \circ^{\zeta^{-3}} }, \qquad \zeta \in \mathbb{G}_9 , \label{hatB3}
\end{align}
where $\mathbb{G}_k$ denotes the set of primitive root of unity of order $k$.

If $l=1$, we have nothing to consider, except that $q^{b_1d_1} \neq 1$; that is, $p^{n}$ does not divide $b_1d_1$.
\medskip

When $l=2$, the braiding is of Cartan type $A_1 \times A_1$, $A_2$, $B_2$, $G_2$, or non-Cartan of type $\hat{B}_2$.
\medskip

$\mathbf{A_1\times A_1}$: We have $1=q_{12}q_{21}=q^{b_1d_2+b_2d_1}$, so $b_1d_2+b_2d_1 \equiv 0 (p^n)$. First of all, $\alpha_1+\gamma_2=\alpha_2+\gamma_1$, because both numbers are less than n; we call m to this number. Also, $a_1c_2+a_2c_1 \equiv 0 (p^{n-m})$. We can describe then the set of solutions as choosing $\alpha_i \leq m <n$, $a_1,a_2,c_1,c_2$ non divisible by $p$ such that $a_1c_2+a_2c_1 \equiv 0 (p^{n-m})$ (one can choose freely three of them and determine the other), and define $\gamma_i=m-\alpha_i$.
\medskip

$\mathbf{A_2}$: Now, $b_1c_1 \equiv b_2c_2 \equiv -b_1c_2-b_2c_1 (p^{n})$, so
  $$ \alpha_1+\gamma_1 = \alpha_2+\gamma_2 = \min \{\alpha_1+\gamma_2, \alpha_2+\gamma_1 \}.$$
From this, $\alpha_1=\alpha_2$ and $\gamma_1=\gamma_2$. Also, $a_1c_1 \equiv a_2c_2 \equiv -a_1c_2-a_2c_1 (p^{n-m})$, where $m= \alpha_i+\gamma_i$. Therefore, $a_1^2+a_1a_2+a_2^2 \equiv 0 (p^{n-m})$, so $p=3$, $m=n-1$ and $a_1 \equiv a_2 (3)$, or $a_1a_2^{-1} \not\equiv 1 (p^{n-m})$ is a cubic root of unity, in which case $p \equiv 1 (3)$.
\medskip

$\mathbf{B_2}$: We have $b_1c_1 \equiv 2b_2c_2 \equiv -b_1c_2-b_2c_1 (p^{n})$, and then $ \alpha_1+\gamma_1 = \alpha_2+\gamma_2 = \min \{\alpha_1+\gamma_2, \alpha_2+\gamma_1 \}$. As before, $\alpha_1=\alpha_2$ and $\gamma_1=\gamma_2$. Also, $a_1c_1 \equiv 2a_2c_2 \equiv -a_1c_2-a_2c_1 (p^{n-m})$, if $m= \alpha_i+\gamma_i$. In this case, $a_1^2+2a_1a_2+2a_2^2 \equiv 0 (p^{n-m})$, so as in \cite{AS1}, this equation has solution if and only if $-1$ is an square modulo $p$, which implies $p \equiv 1 (4)$.
\medskip

$\mathbf{G_2}$: In this case, $b_1c_1 \equiv 3b_2c_2 \equiv -b_1c_2-b_2c_1 (p^{n})$. If $p=3$, then $ \alpha_1+\gamma_1 = \alpha_2+\gamma_2 +1 = \min \{\alpha_1+\gamma_2, \alpha_2+\gamma_1 \}$, which is a contradiction.

If $p\neq 3$, $ \alpha_1+\gamma_1 = \alpha_2+\gamma_2 = \min \{\alpha_1+\gamma_2, \alpha_2+\gamma_1 \}$. Therefore, $\alpha_1=\alpha_2$ and $\gamma_1=\gamma_2$. Also, $a_1c_1 \equiv 3a_2c_2 \equiv -a_1c_2-a_2c_1 (p^{n-m})$ for $m= \alpha_i+\gamma_i$. Therefore, $a_1^2+3a_1a_2+3a_2^2 \equiv 0 (p^{n-m})$, so this equation has solution if and only if $-3$ is an square modulo $p$, which implies $p \equiv 1 (3)$.
\medskip

$\mathbf{\hat{B}_2}$: In this case, $p=3$ and $\zeta$ is a primitive root of order $3^k$ for $2 \leq k \leq n$, so $\xi=\zeta^{\pm 3^{k-1}}$ (note that if $k=1$, then the braiding is of Cartan type $A_2$ or $B_2$). Changing $q$, we can assume $\zeta=q^{3^{n-k}}$, so we have
$$ b_1d_1 \equiv \pm 3^{n-1} (3^n), \quad b_2d_2 \equiv 3^{n-k} (3^n), \quad b_1d_2+b_2d_1 \equiv -3^{n-k} (3^n).  $$
From these equations, $\alpha_1+\gamma_1=n-1, \, \alpha_2+\gamma_2=n-k$ and $\min \{\alpha_1+\gamma_2,\alpha_2+\gamma_1\}=n-k$, so $\alpha_1=\alpha_2$ (in which case $\gamma_1=\gamma_2+k-1$), or $\gamma_1=\gamma_2$ (in which case $\alpha_1=\alpha_2+k-1$), and
$$ a_1c_1 \equiv \pm 1 (3), \quad a_2c_2 \equiv 1 (3^k), \quad \left\{  \begin{array}{ll} a_1c_2+3^{k-1}a_2c_1 \equiv -1 (3^k), & \alpha_1=\alpha_2; \\ 3^{k-1}a_1c_2+a_2c_1 \equiv -1 (3^{k}) , & \gamma_1=\gamma_2. \end{array} \right. $$
Consider the second case; the first is analogue. Multiplying by the invertible element $c_2$ (modulo $p^k$), we have
$$ a_1c_2^2+c_2+3^{k-1}c_1 \equiv 0 (3^k). $$
This equation has a solution if and only if $1-4a_1c_13^{k-1} \equiv 1 \pm 3^{k-1} (3^k)$ is a quadratic residue. Note that
$$ (\pm 3^{k-1}\pm 1)^2 \equiv \pm 2 \cdot 3^{k-1}+1 \equiv 1 \mp 3^{k-1} (3^n),$$
so $1 \pm 3^{k-1}$ are quadratic residues. This provides the possible structures of Yetter-Drinfeld modules of this kind, reconstructing $b_i,d_i$.

\bigskip

Now we are ready to prove the analogous statement to Proposition 5.1 of \cite{AS1} for $p^n$.

\begin{prop}\label{classifnichols}
Let $V$ a Yetter-Drinfeld module over $Z_{p^n}$ of finite Cartan type, $\dim V \geq 3$. Then, $p=3$ and $V$ is of type $A_2\times A_1$ or $A_2\times A_2$.
\end{prop}
\bp
Consider $V$ of dimension 3; we discard first the non-Cartan cases: they are $\hat{B}_2 \times A_1$ or as in \eqref{hatB3}. For all cases, we can consider vertices 1,2 determining a subdiagram of type $\hat{B}_2$, vertices 1,3 not connected, and vertices 2,3 determining a subdiagram of type $A_2$ or $A_1 \times A_1$. From the first condition, $\alpha_1+\gamma_1 > \alpha_2+\gamma_2$; from the second, $\alpha_1+\gamma_1 = \alpha_3+\gamma_3$, and from the last, $\alpha_2+\gamma_2 = \alpha_3+\gamma_3$. But this is a contradiction.

Consider then $V$ of Cartan type. As in \cite{AS1}, it is not of type $A_1\times A_1\times A_1$, so we can assume vertices 1 and 2 of the corresponding Dynkin diagram are connected, and vertices 1 and 3 are disconnected. Moreover, we can assume that if there exists a multiple arrow, it is the one between vertices 1 and 2. That is,
\begin{itemize}
  \item $b_1d_3 \equiv -b_3c_1 (p^n)$;
  \item $b_1c_1 \equiv m b_2c_2 \equiv -b_1c_2-b_2c_1 (p^n)$ for some $m=1,2,3$, in which case $b_2d_3 \equiv -b_3c_2 (p^n)$ (cases $X_2 \times A_1$ for $X=A,B,G$), or $b_3c_3 \equiv b_2c_2 \equiv -b_3c_2-b_2c_3 (p^n)$ and $m=1,2$ (cases $A_3,B_3$), or
  \item $2b_1c_1 \equiv b_2c_2 \equiv -b_1c_2-b_2c_1 (p^n)$, in which case $b_3c_3 \equiv b_2c_2 \equiv -b_3c_2-b_2c_3 (p^n)$ and $m=1,2$ (case $C_3$).
\end{itemize}
As the corresponding submatrices should be of finite Cartan type, we use the previous description for rank $2$. After to reduce the powers of $p$ involved in each equation, we reduce to a equation modulo $p$ for $a_i,c_i$ not divisible by $p$. A detailed study as in \cite{AS1} gives as unique remaining case $A_2\times A_1$, in which case $p=3$.

Thus if we consider $\dim V \geq 4$, each subdiagram of three vertices is of type $A_2\times A_1$, so we have only one possibility: $A_2 \times A_2$ and $p=3$. In this case, we can describe such $V$ as follows:
\begin{align*}
    & b_1=b_2=3^{\alpha}a, \quad d_1=d_2=3^{\gamma}c & \alpha+\gamma=n-1, \, 3 \nmid ac,\\
    & b_3=b_4=3^{\alpha'}a', \quad d_3=d_4=3^{\gamma'}c' & \alpha'+\gamma'=n-1, \, 3 \nmid a'c', \\
    & \alpha+\gamma'=\alpha'+\gamma \Longrightarrow & \alpha=\alpha', \, \gamma=\gamma', \\
    & ac'+a'c \equiv 0 (3).
\end{align*}
\ep

\subsection{Basic quasi-Hopf algebras over $\Z_{p^n}$}

\begin{prop}\label{rank1}
Let $A$ be a basic radically graded Hopf algebra, whit $A[0]=\k[\Z_{p^n}]$ and associator $\Phi_s$, where $p$ does not divide $s$. Then the rank of $A[1]$ over $A[0]$ is $\leq 1$.
\end{prop}
\bp
Suppose there exists $A$ as above such that the rank of $A[1]$ over $A[0]$ is $\geq 2$, and consider $A$ of minimal possible  dimension. By Theorem \ref{thm:projection}, $A=A(H, s)$, $\bar{H}=R\# \Z_{p^{2n}}$ for some Nichols algebra $R$ of diagonal type, $\dim R[1]=2$ and the braiding is given by $(q^{b_id_j})_{i,j=1,2}$. By Heckenberger's classification \cite{H}, it is of Cartan type:
\begin{itemize}
  \item if it of type $A_2$, $B_2$ or $G_2$, then
  \begin{center}
    $b_1d_1+b_1d_2+b_2d_1 \equiv m \, b_2d_2+b_1d_2+b_2d_1 \equiv 0 (p^{2n})$, $m=1,2,3$ respectively;
  \end{center}
  \item if it is of type $A_1 \times A_1$, $b_1d_2+b_2d_1 \equiv 0 (p^{2n})$;
\end{itemize}
or $p=3$ and
\begin{itemize}
  \item it is of standard $\hat{B}_2$ type, with conditions as in Section \ref{classifnichols}.
\end{itemize}
We write $b_i=p^{\alpha_i}a_i$, $d_i=p^{\gamma_i}c_i$, where $p$ does not divide $a_ic_i$. As $b_i \equiv s d_i (p^n)$, we have $\alpha_i=\gamma_i$.

For cases $X_2$, note that $p^{2\alpha_1}a_1c_1 \equiv p^{2\alpha_2}ma_2c_2 (p^{2n})$, so $\alpha_1 =\alpha_2$ (we simply call them $\alpha$). Therefore $a_1c_1 \equiv m a_2c_2 \equiv-a_1c_2-a_2c_1 (p^{2n-2\alpha})$, $a_i \equiv sc_i (p^{n-\alpha})$. These equations imply $a_1^2 \equiv -2a_1a_2 (p^{n-\alpha})$, so $a_1 \equiv -2a_2 (p^{n-\alpha})$, and $ma_2^2 \equiv a_1^2 \equiv 4a_2^2 (p^{n-\alpha})$. That is, $p^{n-\alpha} \mid (4-m)a_2^2$. It follows that $p=3, \alpha=n-1$ and $m=1$. But in this case,
$$ (a_1-a_2)(c_1-c_2) \equiv 3a_1c_1 \, (9) .$$
As $a_1 \equiv -2a_2 \equiv a_2 (9)$ and $a_i \equiv c_i (3)$, it follows that $9 \mid 3a_1c_1$, a contradiction.

For case $A_1 \times A_1$, $\alpha_1 =\alpha_2$ as above, and $a_1c_2+a_2c_1 \equiv 0 (p^{2n-2\alpha})$. It follows that $2a_1a_2 \equiv 0 (p^{n-\alpha})$, which is a contradiction.
\bigskip

From the previous contradictions, the rank of $H[1]$ over $H[0]$ is $\leq 1$.
\ep

\begin{rem}
Note that for any $m$, $A(q)=A(\bB(V)\# \Z_{m^2}, \omega_1)$, where $V$ is the diagonal braided vector space of dimension 1 and braiding given by $q$.
\end{rem}

The question now is what happens when $p$ divides $s$ and we consider the associator given by $\omega_s \in H^3(\Z_m,\k^{\times})$ . Consider the quasi-Hopf algebras $A(H,s)$, for $H=R\# \Z_{p^{2n}}$, and write $s=p^{\theta}t$, where $\alpha \geq 1$ and $p$ does not divide $t$. Consider $\alpha_i, \gamma_i \geq 0$ such that $b_i=p^{\alpha_i}a_i$, $d_i=p^{\gamma_i}c_i$.

\emph{When $R$ has rank one}, the unique condition is $b_1 \equiv sd_1 (p^n)$, which is possible choosing any $d_1$, $\gamma_1 <n$: if $\theta+\gamma_1<n$, then $b_1$ is uniquely defined modulo $p^n$, if $\theta+\gamma_1 \geq n$, simply choose $b_1$ such that $p^n$ divides $b_1$.
\medskip

\emph{When $R$ has rank two}, $R$ is of Cartan type $A_1 \times A_1$, $A_2$, $B_2$ or $G_2$, or $p=3$ and it is of standard type $\hat{B}_2$. In the first case, we will see in Section \ref{nichols} that it is determined by $\alpha_i \leq m <n$, $a_1,a_2,c_1,c_2$ non divisible by $p$ such that $a_1c_2+a_2c_1 \equiv 0 (p^{2n-m})$, and define $\gamma_i=m-\alpha_i$. As also $b_i \equiv sd_i (p^n)$, if we suppose $p^n$ does not divide $b_1$ ($\alpha_1<n$), then $\alpha_1=\theta+\gamma_1$ and $a_1 \equiv c_1 (p)$. But in such case, $\alpha_2= \theta+\gamma_2$ and $a_2 \equiv c_2 (p)$ so $2c_1c_2 \equiv0 (p)$, which is a contradiction. Therefore, $p^n$ divides $b_i$, and we have $\gamma_i+\theta \geq n$. The unique restriction is in consequence to choose $\gamma_i, \alpha_i$ such that $n \leq \gamma_i+\theta$, $\gamma_i+\alpha_i < 2n$.

In the other cases, the condition $b_i \equiv sd_i (p^n)$ gives a contradiction if we suppose $p^n$ does not divide $b_i$ in a similar way to the previous case, so we consider $p^n | b_i$ for all $i$. The other restrictions are given in Section \ref{nichols}. This condition is implicit when we consider the special case for $A_2$ and $p=3$, where $\alpha_i+\gamma_i=2n-1$.
\medskip

\emph{When $R$ has rank greater than 2}, we know $p=3$ and $R$ is of type $A_2\times A_1$ or $A_2\times A_2$. In this case, $\alpha_i+\gamma_i=2n-1$, so any of those examples such that $p^n$ divides $\beta_i$ gives such $H$: we have $\alpha_i=n$, $\gamma_i=n-1$, and the fact that $p$ divides $s$ says that $b_i \equiv sd_i (p^n)$ holds trivially.

\bigskip

From the Theorem \ref{thm:classificationgradedqHA}, we have proved:

\begin{cor}\label{corollary:gradedqHApn}
Let $A=\oplus_{n \geq 0} A[n]$ be a finite dimensional radically graded quasi-Hopf algebra over $\Z_{p^n}$, with associator $\Phi_s$ for some $s$ such that the rank of $A[1]$ over $A[0]$ is $\theta \geq 1$. Then $A= A(\bB(V)\# Z_{p^{2n}},\omega_s)$ for some Yetter-Drinfeld module $V$ over $\Z_{p^{2m}}$ of dimension $\theta=2,3,4$. Moreover, $\theta=3,4$ if and only if $p=3$ and $V$ is of type $A_2 \times A_1$, $A_2 \times A_2$, respectively.
\end{cor}

\bigskip

Also there exist quasi-Hopf algebras $H(p^n,s), \, 1 \leq s\leq p^n-1$, generated by a group-like element $\sigma$ of order $p^n$, with non-trivial associator $\Phi_s$, distinguished elements $\alpha_s=\sigma^{-s},\, \beta=1$ and $S(\sigma)=\sigma^{-1}$. In a similar way as for $n=1$, here an automorphism preserves the power of $p$ which divides s, so we have $2(n-1)$ classes of equivalences up to isomorphism: if $s_0$ is a non quadratic residue coprime with p, then these classes are
$$ H_+(p^n,m):=H(p^n,p^m), \quad H_-(p^n,m):=H(p^n,s_0p^m) \qquad (1 \leq m \leq n-1).$$
Also, these classes are not twist equivalent.

We restrict our attention to the case $p>7$ as above. As a consequence of Theorem \ref{thm:classificationqHA} we have:

\begin{thm}\label{thm:classificationpn}
Let $A$ be a finite dimensional quasi-Hopf algebra such that $$A/Rad(A) \cong \k[\Z_{p^n}],$$ for some prime $p >7$. Then $A$ is twist equivalent to one of the following quasi-Hopf algebras:
\begin{enumerate}
  \item radically graded Hopf algebras $u(\bD,0,0)$ for some datum $\bD$ of type $A_2$, $B_2$ or $G_2$ over $\Z_{p^n}$,
  \item the semisimple quasi-Hopf algebras $H_{\pm}(p^n,m) \, (1 \leq m \leq n-1)$,
  \item the algebras $A(H,\omega_s)$, where $H \cong u(\bD,0,0)$ for some datum $\bD$ of type $A_2$, $B_2$ or $G_2$ over $\Z_{p^{2n}}$ and some $s \in \Upsilon(H)$.
\end{enumerate}
\end{thm}

\end{document}